\documentclass[12pt,reqno,a4paper,twoside]{article}
\usepackage{amsmath,amsthm,amstext,amscd,amssymb,euscript}
\usepackage{epsf}
\renewcommand{\phi}{\varphi}

\renewcommand{\subset}{\subseteq}
\renewcommand{\emptyset}{\varnothing}

\newcommand{\Zd}{\mathbb Z^d}

\newcommand{\dist}{\text{dist}}

\def\1{ {\mathit{1} \!\!\>\!\! I} }

\makeatletter
\@addtoreset{equation}{section}
\makeatother

\newtheorem{ittheorem}{Theorem}
\newtheorem{itlemma}{Lemma}
\newtheorem{itproposition}{Proposition}
\newtheorem{itdefinition}{Definition}
\newtheorem{itremark}{Remark}

\newenvironment{theorem}{\addtocounter{equation}{1}
\begin{ittheorem}}{\end{ittheorem}}

\newenvironment{lemma}{\addtocounter{equation}{1}
\begin{itlemma}}{\end{itlemma}}

\newenvironment{proposition}{\addtocounter{equation}{1}
\begin{itproposition}}{\end{itproposition}}

\newenvironment{definition}{\addtocounter{equation}{1}
\begin{itdefinition}}{\end{itdefinition}}

\newenvironment{remark}{\addtocounter{equation}{1}
\begin{itremark}}{\end{itremark}}

\newcommand{\beq}{\begin{eqnarray}}
\newcommand{\eeq}{\end{eqnarray}}

\newcommand{\be}{\begin{equation}}
\newcommand{\ee}{\end{equation}}

\newcommand{\bl}{\begin{lemma}}
\newcommand{\el}{\end{lemma}}

\newcommand{\br}{\begin{remark}}
\newcommand{\er}{\end{remark}}

\newcommand{\bt}{\begin{theorem}}
\newcommand{\et}{\end{theorem}}

\newcommand{\bd}{\begin{definition}}
\newcommand{\ed}{\end{definition}}

\newcommand{\bp}{\begin{proposition}}
\newcommand{\ep}{\end{proposition}}

\newcommand{\bc}{\begin{corollary}}
\newcommand{\ec}{\end{corollary}}

\newcommand{\bpr}{\begin{proof}}
\newcommand{\epr}{\end{proof}}

\newcommand{\bi}{\begin{itemize}}
\newcommand{\ei}{\end{itemize}}

\newcommand{\ben}{\begin{enumerate}}
\newcommand{\een}{\end{enumerate}}

\newcommand{\Z}{\mathbb Z}
\newcommand{\R}{\mathbb R}
\newcommand{\N}{\mathbb N}

\newcommand{\E}{\mathbb E}

\newcommand{\re}{\ensuremath{\mathcal{R}}}
\newcommand{\ce}{\ensuremath{\mathcal{C}}}

\newcommand{\s}{\ensuremath{\mathcal{S}}}

\newcommand{\fe}{\ensuremath{\mathcal{F}}}

\newcommand{\vi}{\ensuremath{\varphi}}
\newcommand{\La}{\ensuremath{\Lambda}}

\newcommand{\wzero}{{\ensuremath{W \setminus \{ 0 \}}}}
\newcommand{\si}{\ensuremath{\sigma}}
\newcommand{\epsi}{\ensuremath{\epsilon}}

\usepackage{amsmath}
\usepackage{amsfonts}
\usepackage{amsthm}

\newcommand{\eqn}[2]{\begin{equation}\label{#1}#2\end{equation}}
\newcommand{\eqnst}[1]{\begin{equation*}#1\end{equation*}}
\newcommand{\eqnspl}[2]{\begin{equation}\begin{split}\label{#1}%
    #2\end{split}\end{equation}}
\newcommand{\eqnsplst}[1]{\begin{equation*}\begin{split}%
    #1\end{split}\end{equation*}}

\def\conn{\leftrightarrow}
\def\dist{\mathrm{dist}}
\def\diam{\mathrm{diam}}
\def\det{\mathrm{det}}

\def\Prob{\mathrm{Pr}}

\def\LE{\mathrm{LE}}

\def\Zd{{{\mathbb Z}^d}}

\def\Z2{{{\mathbb Z}^2}}

\def\Ed{{{\mathbb E}^d}}

\def\La{\Lambda}

\def\tZd\widetilde{Z^d}

\def\tT{\widetilde{T}}

\def\tW{\widetilde{W}}

\def\hW{\widehat{W}}
\def\tnu{\widetilde{\nu}}

\def\hg{\hat{g}}
\def\eps{\varepsilon}
\def\es{\emptyset}
\def\cnctd{\leftrightarrow}
\def\ntcnctd{\hskip1pt\not\hskip-1pt\leftrightarrow}

\def\cC{{\cal C}}
\def\cD{{\cal D}}
\def\cF{{\cal F}}
\def\cG{{\cal G}}
\def\cH{{\cal H}}
\def\cX{{\cal X}}

\def\cR{{\cal R}}
\def\cQ{{\cal Q}}
\def\cT{{\cal T}}

\def\bv{{\bar v}}

\def\bF{{\bar F}}

\def\hW{{\widehat{W}}}

\def\cT{{\cal T}}
\def\cC{{\cal C}}

\def\CMPsh{Commun.~Math.~Phys.}

\def\JSPsh{J.~Stat.~Phys.}

\def\AOPsh{Ann.~Probab.}

\def\EJPsh{Electron.~J.~Probab.}
\def\MPRFsh{Markov Process. Related Fields}

\begin{document}
\title{{\bf Infinite volume limit
    of the Abelian sandpile model in dimensions $d\geq 3$}}

\author{
Antal A.~J\'arai
\footnote{Carleton University, School of Mathematics and Statistics,
1125 Colonel By Drive, Room 4302 Herzberg Building, Ottawa, ON K1S 5B6,
Canada}\\
Frank Redig
\footnote{Mathematisch Instituut, Universiteit Leiden,
Snellius, Niels Bohrweg 1, 2333 CA Leiden, The Netherlands}\\
}

\maketitle

\footnotesize
\begin{quote}
{\bf Abstract:}
We study the Abelian sandpile model on $\Zd$. In $d\geq 3$ we prove
existence of the infinite volume addition operator, almost surely
with respect to the infinite volume limit $\mu$ of the uniform measures
on recurrent configurations. We prove the existence of a Markov process
with stationary measure $\mu$, and study ergodic properties of this
process. The main techniques we use are a connection between the
statistics of waves and uniform two-component spanning trees and
results on the uniform spanning forest measure on $\Zd$.
\end{quote}
\normalsize

\vspace{12pt}
\section{Introduction}
The Abelian sandpile model (ASM), introduced by Bak, Tang and Wiesenfeld
\cite{btw}, has been studied extensively in the physics literature,
mainly because of its remarkable ``self-organized" critical state.
Many exact results were obtained by Dhar using the group structure of
the addition operators acting on recurrent configurations introduced
in \cite{Dhar}, see for example \cite{Dhar06,Dhar1} for reviews.
The relation between recurrent configurations and spanning
trees, introduced by Majumdar and Dhar \cite{MD92}, has been
used by Priezzhev to compute the stationary height probabilities of the
two-dimensional model in the thermodynamic limit \cite{priezzh94}.
Later on, Ivashkevich, Ktitarev and Priezzhev introduced the concept
of ``waves" to study the avalanche statistics, and made a connection
between two-component spanning trees and waves \cite{IKP94,IKP94b}.
In \cite{priezzh} this connection was used to argue that the critical
dimension of the ASM is $d=4$.

From the mathematical point of view, one is interested in the
thermodynamic limit, both for the stationary measures and for the dynamics.
Recently, in \cite{aj} the connection between recurrent configurations and
spanning trees, combined with results of Pemantle \cite{pmnt91} on
existence and uniqueness of the uniform spanning forest measure on $\Zd$,
has led to the result that the uniform measures $\mu_V$ on recurrent
configurations in finite volume have a unique thermodynamic (weak) limit
$\mu$. (Note: for $d > 4$ the limit was only established for
regular volumes such as cubes centered at the origin. The extension to
arbitrary $V$ is given in the appendix of this paper.)
In \cite{mrs} the existence of a unique limit $\mu$ was proved for an
infinite tree, and a Markov process generated by Poissonian additions
to recurrent configurations was constructed.

A natural continuation of \cite{aj} is therefore to investigate
the dynamics defined on $\mu$-typical configurations. The first question
here is to study the addition operators. We prove that at least in $d\geq 3$,
the addition operators $a_x$, $x \in \Zd$ can be defined on $\mu$-typical
configurations. This turns out to be a rather simple consequence
of the transience of the simple random walk, and we obtain that the
avalanche resulting from adding a particle at a given site will be locally
finite $\mu$-almost surely (all sites topple finitely many times).

Next, in order to construct a stationary process from the infinite
volume addition operators, it is crucial that the measure $\mu$ is
invariant under $a_x$. We show that this is the case if the avalanche
triggered by adding a particle at $x$ is $\mu$-almost surely finite
(only finitely many topplings occur). In order to
establish almost sure finiteness of avalanches, we first prove that the
statistics of waves has a bounded density
with respect to the uniform two-component spanning tree. The final step
then is to show that the component of the uniform two-component spanning
tree corresponding to the wave is almost surely finite in the infinite
volume limit when $d \ge 3$. We deduce this from known results on the
uniform spanning forest \cite{pmnt91,blps01}.
The case $d = 2$ remains an important open
question.

Given existence of $a_x$, and stationarity of $\mu$ under its action,
we can apply the formalism developed in \cite{mrs} to construct a
stationary process which is informally described as follows. Starting
from a $\mu$-typical configuration $\eta$, at each site $x\in\Zd$ grains
are added on the event times of a Poisson process $N^x_t$ with mean
$\vi (x)$, where $\vi (x)$ satisfies the condition
\[
  \sum_{x \in \Zd}\vi (x) G(0,x) <\infty,
\]
with $G$ the Green function of simple random walk in $\Zd$.
The condition ensures that the number
of topplings at $0$ caused by additions at all sites has finite
expectation at any time $t>0$.

In this paper we further study the ergodic properties of the
infinite volume process. We show that tail triviality of the measure
$\mu$ implies ergodicity of the process. We prove that $\mu$ has trivial
tail in any dimension $d \ge 2$. For $2\leq d\leq 4$ this is a rather
straightforward consequence of the fact that the height-configuration
is a (non-local) coding of the edge configuration of the uniform spanning
tree, that is, from the spanning tree in {\em infinite volume} one can
reconstruct the infinite height configuration almost surely.
This is not the case in $d>4$ where we need a separate argument.

Our paper is organized as follows. We start with some notation and
definitions, recalling some basic facts about the ASM. In Sections 3 and 4
we prove existence of the addition operator $a_x$ when $d \ge 3$, and
show invariance of the measure $\mu$, assuming finiteness of avalanches.
In Section 5 we prove existence of inverse addition operators. Sections
6--8 are devoted to establishing finiteness of avalanches in dimensions
$d \ge 3$. In Section 6 we make the precise link between avalanches and
waves, in Section 7 we prove that all waves are finite if the uniform
two-component spanning tree has almost surely a finite component. In Section 8
we prove the required finiteness of the component
in dimensions $d \geq 3$. Finally, in Sections 9 and 10 we discuss tail
triviality of the stationary measure, and
correspondingly, ergodicity of the stationary process.

The review papers \cite{Jarai05,MRS05,Redig06} explain many points
that are presented in less detail here, and may be useful for the reader.

\section{Notation and definitions}
\label{sec:defs}
We consider the Abelian sandpile model, as introduced in \cite{btw}
and generalized by Dhar \cite{Dhar}. One starts from a toppling matrix
$\Delta_{xy}$, indexed by sites in $\Zd$. In this paper $\Delta$ will
always be the degree minus the adjacency matrix (in other words,
minus the discrete lattice Laplacian):
\eqnst
{ \Delta_{xy}
  = \begin{cases}
    2d & \text{if $x=y$,} \\
    -1 & \text{if $|x-y|=1$,} \\
    0  & \text{otherwise.}
  \end{cases} }
A height configuration is a map $\eta: \Zd \to \N = \{1,2,\dots\}$,
and a stable height configuration is such that
$\eta (x) \leq \Delta_{xx}$ for all $x\in\Zd$.
A site where $\eta (x) > \Delta_{xx}$ is called an unstable site.

All stable configurations are collected in the set
$\Omega = \{1, 2, \dots, 2d\}^\Zd$.
We endow $\Omega$ with the product topology. For $V\subset\Zd$,
$\Omega_V = \{1, 2, \dots, 2d\}^V$ denotes the stable configurations
in volume $V$.
If $\eta \in \Omega$ and $W \subset \Zd$, then $\eta_W$ denotes
the restriction of $\eta$ to the subset $W$. We also use $\eta_W$
for the restriction of $\eta \in \Omega_V$ to a subset
$W \subset V$. Given $\eta \in \Omega_V$, $\xi \in \Omega_{V^c}$,
we let $\eta_V \xi_{V^c}$ denote the configuration that agrees
with $\eta$ in $V$ and with $\xi$ in $V^c$.
We define the matrix $\Delta_V$ as the finite volume
analogon of $\Delta$, indexed now by the sites in $V$.
That is, $(\Delta_V)_{xy} = \Delta_{xy}$, $x,y \in V$.
Depending on the context, we sometimes interpret $\Delta_V$ as a
matrix indexed by $\Zd$. In that case
$(\Delta_V)_{xy} = \Delta_{xy} I[x \in V] I[y \in V]$,
where $I[\cdot]$ denotes an indicator function.

The toppling of a site $x$ in finite volume $V$ is defined
on configurations $\eta: V\to \N$ by
\be
T_x (\eta ) (y)  = \eta (y) - (\Delta_V)_{xy}
\ee
A toppling is called legal if the toppled site was unstable,
otherwise it is called illegal.
The stabilization of an unstable configuration
is defined to be the stable result of a sequence
of legal topplings, i.e.,
\be
\s_V (\eta) = T_{x_n}\circ T_{x_{n-1}}\circ\ldots\circ T_{x_1} (\eta),
\ee
where all topplings are legal and such that $\s_V (\eta )$ is
stable. That $\s_V (\eta )$ is well-defined
follows from \cite{Dhar,mrz}, see also \cite{gabr}.
If $\eta$ is stable, then we define $\s_V (\eta ) =\eta$.
The addition operators are then defined by
\be
a_{x,V}\eta = \s_V (\eta + \delta_x ),
\ee
where $\delta_x (y)$ is $1$ for $y = x$ and $0$ otherwise.
As long as we are in finite volume, $a_{x,V}$ is well-defined and
$a_{x,V} a_{y,V} = a_{y,V} a_{x,V}$ (Abelian property).

The dynamics of the finite volume ASM is described as follows:
at each discrete time step choose at random a site $X$
according to a probability measure $p(x)>0, x\in V$, and
apply $a_{X,V}$ to the configuration. After time $n$,
the configuration is $\prod_{i=1}^n a_{X_1,V}\ldots a_{X_n,V} \eta$
where $X_1, \ldots, X_n$ is an i.i.d.~sample from $p$.
This gives a Markov chain with transition operator
\be
Pf(\eta ) = \sum_{x \in V} p(x) f(a_{x,V}\eta )
\ee

A function $f: \Omega \to \R$ is called local, if it only depends
on finitely many coordinates, that is, there exists finite
$V \subset \Zd$, and $g : \Omega_V \to \R$ such that
$f(\eta) = g(\eta_V)$, $\eta \in \Omega$. Local functions are
dense in the space of continuous functions on $\Omega$ with
respect to uniform convergence.

Given a function $F(V)$ defined for all sufficiently large
finite volumes in $\Zd$, and taking values in a metric space
with metric $\rho$, we say that $\lim_V F(V) = a$, if for all
$\eps > 0$ there exists finite $W$, such that
$\rho(F(V),a) < \eps$ whenever $V \supseteq W$.
For a probability measure $\nu$, $\E_\nu$ will denote
expectation with respect to $\nu$. The boundary of $V$ is defined by
$\partial V = \{y \in V : \text{$y$ has a neighbour in $V^c$} \}$,
while its exterior boundary is defined by
$\partial_e V = \{y \in V^c : \text{$y$ has a neighbour in $V$} \}$.

\subsection{Recurrent configurations}
A stable configuration $\eta \in \Omega_V$ is called
recurrent, if it is recurrent in the Markov chain defined in Section
\ref{sec:defs}. Equivalently, $\eta$ is recurrent, if for any $x \in V$
there exists $n = n_{x,\eta}$ such that $a_{x,V}^n \eta = \eta$.
We denote by $\re_V$ the set of recurrent configurations.
The addition operators $a_{x,V}$ restricted to $\re_V$ have
well-defined inverses $a_{x,V}^{-1}$, and therefore
form an Abelian group under composition. From this fact one easily
concludes that the uniform measure $\mu_V$ on $\re_V$ is the unique
invariant measure of the Markov chain.
One can compute the number of recurrent configurations:
\be\label{Dharformula1}
  |\re_V| = \det (\Delta_V),
\ee
see \cite{Dhar}.
Another important identity of \cite{Dhar} is the following.
Denote by $N_V (x,y,\eta)$ the number of legal topplings at
$y$ needed to obtain $a_x\eta$ from $\eta + \delta_x$. Then
the expectation of $N_V$ satisfies
\be\label{klingsor}
  \E_{\mu_V} ( N_V(x,y,\eta) )
  = G_V (x,y)
  \stackrel{\mathrm{def}}{=} (\Delta_V^{-1})_{xy}
  \qquad \text{(Dhar's formula).}
\ee
From this and the Markov inequality, one also obtains $G_V (x,y)$
as an estimate of the $\mu_V$-probability
that a site $y$ has to be toppled if one adds at $x$.
We also note that for our specific choice of $\Delta$, $G_V$ is
$(2d)^{-1}$ times the Green function of simple random walk in
$V$ killed upon exiting $V$.

Recurrent configurations are characterized
by the so-called burning algorithm \cite{Dhar}.
A configuration $\eta$ is recurrent
if and only if it does not contain a so-called
forbidden sub-configuration, that is, a subset
$W \subset V$ such that for all $x \in W$:
\be
\eta (x) \leq - \sum_{y \in W \setminus \{x\}} \Delta_{xy}.
\ee

From this explicit characterization, one easily infers
a consistency property:
if $\eta\in\re_V$ and $W\subset V$ then
$\eta_W\in\re_W$. This suggests a natural definition of
``recurrent configurations in infinite volume'': we
say that $\eta \in \Omega$ is recurrent, if its restriction
to any finite $V$ is. We denote this set by $\re$:
\eqnst
{ \re
  = \{ \eta \in \Omega : \text{$\eta_V \in \re_V$ for
    all finite $V \subset \Zd$} \}. }

\subsection{Infinite volume: basic questions and results}
In studying infinite volume limits of the ASM, the following
questions can be addressed. In this (non-exhaustive) list, any
question can be asked only after a positive answer to all
previous questions.
\ben
\item Do the measures $\mu_V$ weakly converge to a measure
  $\mu$? Does $\mu$ concentrate on the set $\re$?
\item Is the addition operator $a_x$ defined
  on $\mu$-a.e.~configuration $\eta\in\re$, and
  does it leave $\mu$ invariant?
  Does the Abelian property still hold in infinite volume?
\item Can one define a natural Markov process on $\re$ with stationary
  distribution $\mu$?
\item Does the stationary Markov process of question 3 have
  good ergodic properties?
\een

Question 1 is easily solved for the one-dimensional lattice
${\mathbb Z}$, however, $\mu$ is trivial, concentrating on the single
configuration that is identically $2$.
Hence no further questions on our list are relevant in that case.
See \cite{MRSvM} for a result on convergence to equilibrium
in this case. For an infinite regular tree, the first three
questions have been answered affirmatively and the fourth question
remained open \cite{mrs}. For dissipative models, that is when
$\Delta_{xx} > 2d$, all four questions are affirmatively answered
when $\Delta_{xx}$ is sufficiently large \cite{mrs2}.

For $\Zd$, question 1 is positively answered in any dimension
$d \ge 2$, using a correspondence between spanning trees and
recurrent configurations and properties of the uniform
spanning forest on $\Zd$ \cite{aj}. The limiting measure
$\mu$ is translation invariant. The proof of convergence in
\cite{aj} in the case $d > 4$ is restricted to regular
volumes, such as a sequence of cubes centered at the origin.
In the appendix, we prove convergence along an
arbitrary sequence of volumes using a random walk result \cite{jk}.

In this paper we study questions 2, 3 and 4 for
$\Zd$, $d\geq 3$, and all questions are affirmatively
answered.

The main problem is to prove that avalanches are
almost surely finite. This is done by a decomposition of avalanches
into a sequence of waves (cf.~\cite{IKP94b,priz}), and studying
the almost sure finiteness of the waves.
The latter can be achieved by a two-component
spanning tree representation of waves, as
introduced in \cite{IKP94b,priz}.
We then study the uniform two-component spanning
tree in infinite volume and prove that
the component containing the origin is
almost surely finite. This turns out to be sufficient
to ensure finiteness of all waves.

\section{Existence of the addition operator}
\label{exax}
In this section we show convergence of
the finite volume addition operators
to an infinite volume addition operator when $d \geq 3$.
This turns out to be easy, but in order to make
appropriate use of this infinite volume addition
operator, we need to establish that $\mu$ is invariant
under its action, and for the latter we need to show that
avalanches are finite $\mu$-a.s.

Let $a_{x,V}$ denote the addition operator acting on
$\Omega_V$. We define a corresponding operator acting
on $\Omega$ using the finite $(V)$-volume rule,
that is, grains falling out of $V$ disappear.
More precisely, for $\eta \in \Omega$ and $V\ni x$, we define (with
slight abuse of notation)
\be
\label{axv}
  a_{x,V} \eta = (a_{x,V} \eta_V) \eta_{V^c}.
\ee
Given $\eta \in \Omega$, call $N_V(x,y,\eta)$ the
number of topplings caused at $y$ by addition at
$x$ in $\eta$, using the finite $(V)$-volume rule.
Then
\be
  \eta + \delta_x - \Delta_V N_V(x,\cdot,\eta)
  = a_{x,V} \eta, \qquad
  \eta \in \Omega, x \in V
\ee
where $\Delta_V$ is indexed by $\Zd$.

We start with the following simple lemma:
\bl
\label{simple}
$N_V(x,y,\eta)$ is a non-decreasing function of
$V$ and depends on $\eta$ only through $\eta_V$.
\el

\bpr
Let $V\subset W$. Suppose we add a grain
at $x$ in configuration $\eta$. We perform
topplings inside $V$ until inside $V$ the
configuration is stable, using the finite
$(W)$-volume rule. The result of this procedure is
a configuration $(a_{x,V}\eta_V)\xi_{V^c}$,
where possibly $\xi_{V^c\cap W}$ is not stable.
In that case we perform all the necessary topplings still
needed to stabilize $(a_{x,V}\eta_V)\xi_{V^c\cap W}$ inside
$W$, using the finite $(W)$-volume rule. This can only cause
potential extra topplings at any site $y$ inside $V$.
\epr

From Lemma \ref{simple} and by monotone convergence:
\be
\E_\mu (\sup_V N_V (x,y,\eta)) = \lim_{V} \E_\mu (N_V (x,y,\eta)).
\ee
By weak convergence of $\mu_V$ to $\mu$, and by Dhar's formula \eqref{klingsor}:
\eqnspl{e:bound-by-G}
{ \lim_{V} \E_\mu (N_V (x,y,\eta))
  &=
  \lim_{V} \lim_{W\supseteq V} \E_{\mu_W} (N_V (x,y,\eta)) \\
  &\leq
  \lim_{V} \lim_{W\supseteq V} \E_{\mu_W} (N_W (x,y,\eta)) \\
  &=
  \lim_{W} G_W (x,y) = G(x,y), }
where $G(x,y) = \Delta^{-1}_{xy}$ is $(2d)^{-1}$ times
the Green function of simple random walk in $\Zd$.
In the last step we used that $d \geq 3$, otherwise
$G_W (x,y)$ diverges as $W\uparrow \Zd$.
This proves that for all $x,y\in\Zd$,
$N(x,y,\eta)=\sup_V N_V (x,y,\eta)$ is $\mu$-a.s.~finite.
Hence
\be
\label{good-etas}
\mu \left( \forall x,y\in \Zd: N(x,y,\eta) < \infty \right) = 1.
\ee
Therefore, on the event in \eqref{good-etas}, we can define
\be\label{axdef}
a_x \eta = \eta + \delta_x - \Delta N(x,\cdot,\eta).
\ee
It is easy to see that $a_x \eta$ is stable. This is because
$a_x \eta (y)$ is already determined by the number of topplings
at $y$ and its neighbours, and this is the same as it was in
a large enough finite volume $V$. By similar reasons, we also get
\be
\label{limaxv}
  a_x \eta
  = \lim_{V} a_{x,V} \eta, \quad\text{$\mu$-a.s.,}
\ee
where $a_{x,V}$ is defined in \eqref{axv}.

From its definition, one sees that $a_x$ is well behaved
with respect to translations. Let
$\theta_x$ denote the shift on configuration, that is,
$\theta_x\eta (y) =\eta (y-x)$. Then
\be
\label{cov}
a_{x} = \theta_{x} \circ a_0 \circ \theta_{-x},
\ee
whenever either side is defined.

Note that with the above definition of $a_x$,
there can be infinite avalanches.
However, if the volume increases, it cannot happen that the number
of topplings at a fixed site diverges, and that is the only
problem for defining $a_x$ (a problem which may arise in $d = 2$).
More precisely, an infinite avalanche that
leaves eventually every finite region does not pose a problem
for defining the addition operator. However, as we will see
later on, infinite avalanches do cause problems in defining
a stationary process, at least with our current methods.
Hence we need extra arguments to show that the total number
of topplings is finite $\mu$-a.s.

We define the avalanche cluster caused by addition at
$x$ to be the set
\be
  \ce_x(\eta ) = \{ y\in \Zd: N(x,y,\eta) > 0 \}
\ee
We say that the avalanche at $x$ is finite in $\eta$ if
$\ce_x(\eta )$ is a finite set.
We say that $\mu$ has the finite avalanche property,
if for all $x\in\Zd$, $\mu (|\ce_x| < \infty) = 1$.

In Sections 6--8, we prove the following theorem:

\bt
\label{thm:finite-aval}
Assume $d \geq 3$. Then $\mu$ has the finite avalanche
property, that is, $\mu( |C_x| < \infty ) = 1$ for all
$x \in \Zd$.
\et

In Section 4, we show that Theorem \ref{thm:finite-aval}
has the following consequence.

\bp
\label{prop:inv}
Assume $d \geq 3$. Then $\mu$ is invariant under the
action of $a_x$, $x \in \Zd$, that is, for any
$\mu$-integrable function $f$ and for any $x \in \Zd$,
\be
  \int f(a_x\eta ) d\mu = \int f(\eta) d\mu.
\ee
\ep

Before moving on to the proofs of Theorem \ref{thm:finite-aval}
and Proposition \ref{prop:inv}, we prove some of their
easy consequences.

Integrating \eqref{axdef} over $\mu$ and using
Proposition \ref{prop:inv}, we easily obtain the
following infinite volume analogue of Dhar's formula.

\bp Assume $d \geq 3$. Then
\be
\E_\mu (N(x,y,\eta)) = G(x,y)
\ee
\ep

At this point, we cannot compose different
$a_x$, since $a_x$ is only defined almost surely.

\bp
\label{prop:compose}
Assume $d \geq 3$. There exists a $\mu$-measure one set
$\Omega' \subset \re$ with the following properties.
\begin{itemize}
\item[(i)] For any $\eta \in \Omega'$ and $x \in \Zd$, there
  exists finite $V_x(\eta) \subset \Zd$, such that for
  all $W \supseteq V_x (\eta )$
\eqnst
{ a_x \eta = a_{x,W}\eta. }
\item[(ii)] For any $\eta \in \Omega'$ and $x \in \Zd$ we
  have $a_x \eta \in \Omega'$.
\end{itemize}
Consequently, for any $\eta\in\Omega'$, and any
$x_1,\ldots,x_n\in\Zd$, $a_{x_n} a_{x_{n-1}} \dots a_{x_1} \eta$
is well-defined and all avalanches involved are finite.
\ep

\bpr
Define
\eqnst
{ \Omega_0
  = \{ \eta \in \re: \text{$|\ce_x(\eta)| < \infty$ for all
    $x \in \Zd$} \}. }
By Theorem \ref{thm:finite-aval} and since $\mu$ is concentrated
on $\re$, we have $\mu(\Omega_0) = 1$.
Property (i) in the proposition is satisfied for all
$\eta \in \Omega_0$. For (ii), we need to find a subset of $\Omega_0$
invariant under all the $a_x$'s. For $n \ge 1$, define
inductively the sets
\eqnst
{ \Omega_n
  = \Omega_{n-1} \cap
    \bigcap_{x \in \Zd} a_x^{-1} ( \Omega_{n-1} ), }
where $a_x^{-1}$ here denotes inverse image (not to be confused
with the inverse operator defined later). Since the
$a_x$ are measure preserving, it follows by induction that
$\mu(\Omega_n) = 1$ for all $n$. Also, $a_x$ maps $\Omega_n$
into $\Omega_{n-1}$. Therefore,
$\Omega' = \cap_{n \ge 0} \Omega_n$ satisfies both properties
stated.
\epr

The following proposition shows that
the Abelian property holds in infinite volume.

\bp
Assume $d \geq 3$. Then
\be
  a_x a_y \eta = a_y a_x \eta, \qquad \eta \in \Omega'.
\ee
\ep

\bpr
By Proposition \ref{prop:compose}, for $\eta\in\Omega'$
and for $W \supseteq V_y(\eta) \cup V_x(a_y \eta)
\cup V_x(\eta) \cup V_y (a_x\eta)$ we have
\be
a_x a_y \eta = a_{x,W} a_{y,W} \eta = a_{y,W} a_{x,W} \eta
= a_y a_x \eta.
\ee
\epr

\section{Invariance of $\mu$ under $a_x$}
\label{axinv}
In this section, we show that $\mu$ is invariant
under the addition operators, if it has the
finite avalanche property, that is, we show that
Theorem \ref{thm:finite-aval} implies
Proposition \ref{prop:inv}.

\emph{Proof of Proposition \ref{prop:inv} (assuming
Theorem \ref{thm:finite-aval}).}
It is enough to prove the claim for $f$ a local function.
In that case, we have
\beq
  \int f( a_x \eta ) d\mu
  &= & \int f(a_{x,V}\eta )d\mu + \epsilon_1 (V,f)
    \nonumber\\
  &= & \int f(a_{x,V}\eta )d\mu_W + \epsilon_1 (V,f)
    + \epsilon_2 (V,W,f) \nonumber\\
  &= & \int f(a_{x,W}\eta )d\mu_W + \epsilon_1 (V,f)
    + \epsilon_2 (V,W,f) + \epsilon_3 (V,W,f) \nonumber
\eeq
Here $\epsi_1$ and $\epsi_2$ can be made arbitrarily small
by \eqref{limaxv} and by weak convergence. We also have
\be
  |\epsi_3(V,W,f)|
  \le 2 \|f\|_\infty \mu_W (a_{x,W} f\not= a_{x,V} f). \nonumber
\ee
Next, by invariance of $\mu_W$ under the action of $a_{x,W}$,
\be
  \int f(a_{x,W}\eta ) d\mu_W
  = \int f d\mu_W
  = \int f d\mu + \epsi_4(W,f).
\ee
Here, by weak convergence, $\epsi_4$ can be made arbitrarily small.
Therefore, combining the estimates, we conclude
\be\label{kwakk}
  \left| \int f(a_x\eta ) d\mu - \int f(\eta ) d\mu \right|
  \leq C \limsup_{V} \limsup_{W \supseteq V}
    \mu_W ( a_{x,W} f\not= a_{x,V} f ).
\ee

Define the avalanche cluster in volume $W$ by
\eqnst
{ \ce_{x,W}(\eta) = \{ y \in W : N_W(x,y,\eta) > 0 \},
    \quad \eta \in \Omega. }
Let $D_f$ denote the dependence set of the local function $f$.
On the event $\ce_{x,W}(\eta) \cap \partial V = \es$ we have
$a_{x,V} \eta = a_{x,W} \eta$. Hence, provided $D_f \subset V$,
we have
\eqnst
{ \mu_W (a_{x,W} f\not= a_{x,V} f)
  \le \mu_W ( \ce_{x,W} \cap \partial V \not= \es ). }
The event on the right hand side is a cylinder event (only depends
on heights in $V$). Therefore, the right hand side
approaches $\mu ( \ce_x \cap \partial V \not= \es )$, as
$W \uparrow \Zd$. By Theorem \ref{thm:finite-aval},
\eqnst
{ \lim_V \mu ( \ce_x \cap \partial V \not= \es )
  = \mu ( |\ce_x| = \infty ) = 0, }
which completes the proof.

\section{Inverse addition operators}

In this section we prove that $a_x$ has an inverse defined
$\mu$-a.s., provided $\mu$ has the finite avalanche property.
In fact, we show that $a_x^{-1}$ is the limit of finite volume
inverses.
Let $a_{x,V}^{-1}$ denote the inverse of $a_{x,V}$ acting
on $\re_V$. We define a corresponding operator acting on
$\re$, by
\[
  a_{x,V}^{-1}\eta = (a_{x,V}^{-1}\eta_V)\eta_{V^c},
  \quad \eta \in \re.
\]
This is well-defined, since $\eta_V\in\re_V$.

Recall that under the finite avalanche property,
Proposition \ref{prop:compose} provided a set $\Omega'$ of
recurrent configurations such that for any $\eta \in \Omega'$
and every $x\in\Zd$, there exists a finite set $V_x(\eta)$
such that $a_x\eta = a_{x,V_x(\eta)}\eta$.

\bp
Assume $d \geq 3$. There exists a $\mu$-measure one set
$\Omega'' \subset \Omega'$ with the following properties.
\begin{itemize}
\item[(i)] For any $\eta \in \Omega''$ and $x \in \Zd$ there
  exists finite ${\bar V} = {\bar V}_x (\eta)$ such that
  $a_{x,W}^{-1}\eta = a_{x,{\bar V}}^{-1} \eta$
  for all $W \supseteq {\bar V}$.
\item[(ii)] If we define
  $a_x^{-1} \eta = a_{x,{\bar V}_x(\eta)}^{-1} \eta$, then
  $a^{-1}_x a_x \eta = a_x a_x^{-1} \eta = \eta$ for
  $\eta \in \Omega''$.
\item[(iii)] As operators on $L_2 (\mu)$, $a_x^*= a_x^{-1}$,
  that is, the $a_x$ are unitary operators.
\end{itemize}
\ep

\bpr
We first prove that
\be
\label{invstable}
  \lim_{{\bar V}} \mu \left( \exists W \supseteq {\bar V} :
    a_{x,W}^{-1} \eta \not= a_{x,{\bar V}}^{-1} \eta \right)
  = 0.
\ee
We have
\eqnspl{invcomp}
{ &\mu \left( \exists W \supseteq {\bar V} :
    a_{x,W}^{-1} \eta \not= a_{x,{\bar V}}^{-1} \eta \right) \\
  &\qquad = \mu \left( \exists W \supseteq {\bar V} :
    a_{x,W}^{-1} a_x \eta \not= a_{x,{\bar V}}^{-1} a_x \eta \right) \\
  &\qquad = \mu \left( \exists W \supseteq {\bar V} :
    \text{$a_{x,W}^{-1} a_x \eta \not= a_{x,{\bar V}}^{-1} a_x \eta$
    and $\forall W \supseteq {\bar V}$,
    $a_{x,W} \eta = a_{x,{\bar V}} \eta$} \right) \\
  &\qquad\qquad + \epsi_{{\bar V}} \\
  &\qquad = \epsi_{{\bar V}}. }
Here we used the invariance of $\mu$ under $a_x$ in the first step.
The last step follows because if
$a_x \eta = a_{x,W} \eta = a_{x,{\bar V}} \eta$,
then
\be
\label{eq:inv-rel}
  a_{x,W}^{-1} a_x \eta
  = a_{x,W}^{-1} a_{x,W} \eta
  = \eta
  = a_{x,{\bar V}}^{-1} a_{x,{\bar V}} \eta
  = a_{x,{\bar V}}^{-1} a_x \eta.
\ee
As for $\epsi_{{\bar V}}$ we have
\be
  \epsi_{{\bar V}}
  \leq \mu \left ( \exists W \supseteq {\bar V} :
    a_{x,W} \eta \not= a_{x,{\bar V}} \eta \right)
\ee
which converges to zero as ${\bar V} \uparrow\Zd$,
by the finite avalanche property.

Since the events in \eqref{invstable} are decreasing
in ${\bar V}$, we get that property (i) of the
proposition is satisfied for $\mu$-a.e.~$\eta$.
Let us define $a_x^{-1}$ on this set by setting
$a_x^{-1} \eta = a_{x,{\bar V}_x(\eta)}^{-1} \eta$.
Then \eqref{eq:inv-rel} shows that $a_x$
has an inverse on a full measure set. By standard
arguments, similar to the one in Proposition
\ref{prop:compose}, we can shrink the set $\Omega'$
appropriately to a set $\Omega''$ of full measure
such that $\Omega''$ is invariant under
$a_x$ and $a_x^{-1}$ for all $x \in \Zd$. Then
(i) and (ii) will hold for $\Omega''$.

The last statement of the proposition is an obvious consequence
of the first two.
\epr

The above proposition proves that as operators
on $L_2 (\mu )$, the $a_x$ generate an (Abelian) unitary group,
which we denote by $G$.

\section{Waves and avalanches}
\label{sec:waves}
The goal of Sections \ref{sec:waves}, \ref{sec:finite} and
\ref{sec:finiteness} is to prove Theorem \ref{thm:finite-aval},
saying that $\mu$ has the finite avalanche property.

We will decompose avalanches into so-called waves, that
correspond to carrying out topplings in a special order.
We prove that almost surely, there is a finite number
of waves, and that all waves are almost
surely finite. Without loss of generality, we assume that
the particle is added at $x = 0$, the origin. Since the
site where we add will remain fixed throughout
Sections \ref{sec:waves}, \ref{sec:finite} and
\ref{sec:finiteness}, henceforth we drop
indices referring to $x$ from our notation, and simply
write $\ce = \ce(\eta)$ for the avalanche cluster
$\ce_0(\eta)$.

We start by recalling the definition of a wave
from \cite{IKP94b,priz}. Consider a finite volume $W \ni 0$,
and add a grain at site $0$ in a stable configuration.
If the site becomes unstable, then topple it once
and topple all other sites that become unstable, except $0$.
It is easy to see that in this procedure a site can topple
at most once. The set of toppled sites is called the first wave.
Next, if $0$ has to be toppled again, we start a second wave,
that is we topple $0$ once again, and topple all other sites
that become unstable. We continue as long as $0$ needs to
be toppled.

We define $\alpha_{W} (\eta)$ to be
the number of waves caused by addition at $0$
in the volume $W$.
By definition, $\alpha_{W}$ is the number
of topplings at $0$,
that is $\alpha_{W}(\eta) = N_W(0,0,\eta)$.
For fixed $W$, let $\ce_{W}(\eta)$ denote the
avalanche cluster in volume $W$. We decompose
$\ce_{W}$ as
\be
\label{e:wavedecomp}
  \ce_{W} (\eta )
  = \bigcup_{i=1}^{\alpha_{W} (\eta)} \Xi_{W}^i (\eta),
\ee
where $\Xi_{W}^i (\eta)$ is the $i$-th wave in $W$ caused
by addition at $0$.

We can define waves in infinite volume as we defined
the toppling numbers and avalanches in
Section \ref{exax}, by monotonicity in the
volume. More precisely, the definition is as follows.
Fix $\eta \in \Omega$, and assume that $0$ is
unstable in $\eta + \delta_0$.
By the argument of Lemma \ref{simple} in Section \ref{exax},
$\Xi_{W}^1(\eta)$ is non-decreasing in $W$, and therefore we can define
$\Xi^1(\eta) = \cup_W \Xi_{W}^1(\eta)$.
Let $\eta_1$ denote the configuration obtained by toppling
every site in $\Xi^1(\eta)$ once, that is
$\eta_1 = \lim_W \left[ \prod_{x \in \Xi^1_W(\eta)} T_x \right] \eta$.
Note that carrying out the first wave in the unstable configuration
$\eta + \delta_0$ results in exactly $\eta_1 + \delta_0$.
By the definition of a wave, all sites are stable in $\eta_1$.
If $0$ is unstable in $\eta_1 + \delta_0$, we consider
$\Xi^1_W(\eta_1)$.
Since this is again nondecreasing in $W$, we can define
$\Xi^2(\eta) = \cup_W \Xi^1_W (\eta_1)$.
This is the second wave in infinite volume.
Let $\eta_2 = \lim_W \left[ \prod_{x \in \Xi^1_W(\eta_1)} T_x \right] \eta_1$.
Then $\eta_2 + \delta_0$ is the result of carrying out the
first two waves in infinite volume on $\eta + \delta_0$.
Note that if $|\Xi^1(\eta)| < \infty$, we have
$\Xi^1_W (\eta_1) = \Xi^2_W (\eta)$
for all large $W$, and consequently,
$\Xi^2(\eta) = \lim_W \Xi^2_W(\eta)$.
We similarly define inductively
$\Xi^i(\eta) = \cup_W \Xi^1_W (\eta_{i-1})$ and
$\eta_i = \lim_W \left[ \prod_{x \in \Xi^1_W(\eta_{i-1})} T_x \right] \eta_{i-1}$.
Under the assumption $|\Xi^j(\eta)| < \infty$, $1 \le j < i$, we also have
$\Xi^i(\eta) = \lim_W \Xi^i_W (\eta)$. For convenience, we define
$\Xi^i_W$ or $\Xi^i$ as the empty set, whenever such waves
do not exist.

The easy part in proving finiteness of
avalanches is to show that the number of waves
is finite. Since $\alpha_{W}(\eta)$ is non-decreasing
in $W$, it has a pointwise limit $\alpha (\eta)$,
and as in \eqref{e:bound-by-G},
\be
  \E_\mu (\alpha) \leq  G(0,0) < \infty.
\ee
This implies $\alpha < \infty$ $\mu$-a.s.

In order to prove that $\ce (\eta)$ is finite
$\mu$-a.s., we show, by induction on $i$, that
all sets $\Xi^i (\eta)$ are finite $\mu$-a.s.
We base the proof on the following proposition, proved in
Sections \ref{sec:finite} and \ref{sec:finiteness}.
\bp
\label{prop:finite}
Assume $d \geq 3$. For $i \ge 1$ we have
\be
\label{lap}
  \lim_{V} \limsup_{W \supseteq V}
     {\mu_W} ( \Xi_{W}^i \not\subset V)
  = 0.
\ee
\ep

\emph{Proof of Theorem \ref{thm:finite-aval} (assuming
Proposition \ref{prop:finite}).}
We prove by induction on $i$ that $\mu(|\Xi^i| < \infty) = 1$,
$i \ge 1$.
To start the induction, note that
$\{ \Xi^1 \subset V \}$ is a local event, hence
weak convergence of $\mu_W$ to $\mu$ and
Proposition \ref{prop:finite} with $i=1$ imply that
$\mu(|\Xi^1| < \infty) = 1$.
Assume now that $\mu(|\Xi^j| < \infty) = 1$, $1 \le j < i$.
Then
\be
\label{e:iandj}
  \mu( \Xi^i \not\subset V )
  \le \mu( \Xi^i \not\subset V,\, \Xi^j \subset V',\, 1 \le j < i )
    + \mu(\Xi^j \not\subset V' \text{ for some $1 \le j < i$}).
\ee
By the induction hypothesis, the second term on the right hand side
can be made arbitrarily small by choosing $V'$ large. For fixed
$V'$, the event in the first term is a local event (only depends on
sites in $V' \cup \partial_e V' \cup V \cup \partial_e V$). Therefore,
the first term on the right hand side of \eqref{e:iandj} equals
\be
  \lim_W \mu_W ( \Xi_{W}^i \not\subset V,\,
    \Xi_{W}^j \subset V',\, 1 \le j < i )
  \le \limsup_W \mu_W ( \Xi_{W}^i \not\subset V ).
\ee
Here the right hand side goes to $0$ as $V \uparrow \Zd$, by
Proposition \ref{prop:finite}, showing that
$\mu(|\Xi^i| < \infty) = 1$.

On the event
$\{ \alpha < \infty \} \cap \{ |\Xi^i| < \infty,\, i \ge 1 \}$,
we can pass to the limit in \eqref{e:wavedecomp} and obtain the
decomposition
\be
  \ce (\eta)
  = \bigcup_{i=1}^{\alpha (\eta)} \Xi^i (\eta ).
\ee
It follows that $\mu( |\ce| < \infty ) = 1$, which completes
the proof of Theorem \ref{thm:finite-aval}.

\section{Finiteness of waves}
\label{sec:finite}
In this section we prove Proposition \ref{prop:finite} saying
that waves are finite. The proof is based on a correspondence with
two-component spanning trees due to
Ivashkevich, Ktitarev and Priezzhev \cite{IKP94b}, which is
recalled below. The correspondence
allows us to use some results on the uniform spanning forest
that are stated separately as Proposition \ref{prop:weaklimit}
below. The argument is completed with a proof of
Proposition \ref{prop:weaklimit} in Section \ref{sec:finiteness}.

We now describe the representation of waves as two-component
spanning trees from \cite{IKP94,IKP94b}.
Consider a configuration $\eta_W \in \re_W$ with
$\eta_W(0) = 2d$, and suppose we add a particle at $0$.
Consider the first wave, which is entirely determined by the
recurrent configuration $\eta_\wzero \in \re_\wzero$.
The result of the first wave on $\wzero$ is given by
\be
\label{e:1stwave}
  S^1_{W} (\eta)
  = \Big( \prod_{j \sim 0} a_{j,\wzero} \Big) \eta_{\wzero}
  \in \re_\wzero.
\ee
We associate to any $\xi \in \re_\wzero$ a tree
$\cT_{W}(\xi)$. The tree will represent a wave starting at $0$
in $\xi$. For the definition of the tree, we use Majumdar and
Dhar's tree construction \cite{MD92}.

Denote by $\hW$ the graph obtained from $\Zd$ by identifying
all sites in $\Zd \setminus (\wzero)$ to a single site $\delta_\hW$
(removing loops). By \cite{MD92}, there is a one-to-one map between
recurrent configurations $\xi \in \re_\wzero$ and spanning trees
of $\hW$. The correspondence is given by following the spread of an
avalanche started at $\delta_\hW$. Initially, each
neighbour of $\delta_\hW$ receives a number of grains equal to the
number of edges connecting it to $\delta_\hW$, which results in every
site in $W \setminus \{ 0 \}$ toppling exactly once.
The spanning tree records the sequence
in which topplings have occurred. There is some flexibility in how to
carry out the topplings (and hence in the correspondence with spanning
trees), and here we make a specific choice in accordance with
\cite{IKP94b}. Namely, we first transfer grains from $\delta_\hW$ only
to the neighbours of $0$, and carry out all possible topplings.
We call this the first phase. When we apply the process to
$\xi = \eta_\wzero$, the set of sites that topple in the first
phase is precisely $\Xi^1_W(\eta) \setminus \{ 0 \}$.
Next we transfer grains from $\delta_\hW$ to the boundary sites of $W$,
which will cause topplings at all sites that were not in the wave;
this is the second phase.

The two phases can alternatively be described via the burning
algorithm of Dhar \cite{Dhar}, which in the above
context looks as follows. For convenience, let $\tW$ denote the
graph obtained by identifying all sites in $\Zd \setminus W$ to
a single site $\delta_W$. That is, $\hW$ can be obtained from $\tW$
by identifying $0$ and $\delta_W$, and calling it $\delta_\hW$.
We start with all sites of $\tW$ declared unburnt.
At step $0$ we burn $0$ (the origin). At step $t$, we
\be
\label{e:step}
  \text{burn all sites $y$ for which $\xi(y) > $ current
   number of unburnt neighbours of $y$.}
\ee
The process stops at some step $T = T(\xi)$.
The sites that burn up to time $T$ is precisely the
sites toppling in the first phase. We continue by burning $\delta_W$
in step $T+1$, and then repeating \eqref{e:step} as long as there
are unburnt sites.

Following Majumdar and Dhar's construction \cite{MD92}, we
connect with an edge each $y \in \wzero$ burnt at time $t$ to
a unique neighbour $y'$ (called the parent of $y$) burnt at time
$t-1$. This defines a spanning subgraph of $\tW$ with two
tree components, having roots $0$ and $\delta_W$.
We denote by $\cT_{W}(\xi)$ the tree component having root $0$.
Since $\cT_W(\xi)$ does not contain the vertex $\delta_W$,
it can be identified with a subgraph of $\Zd$, and we will
do so in what follows.
(Identifying $0$ and $\delta_W$ merges the two trees
into a spanning tree of $\hW$, yielding the usual spanning tree
representation of $\xi$.)

With slight abuse of language, we refer to the
two-component spanning subgraph as a two-component spanning tree.
By observations made earlier, when $\xi = \eta_\wzero$,
the vertex set of $\cT_W(\eta_\wzero)$ is the first wave
$\Xi^1_W(\eta)$.

We can generalize the above construction to further waves as follows.
We define
\be
\label{e:kthwave}
  S^k_{W} (\eta)
  = \Big( \prod_{j \sim 0} a_{j,\wzero} \Big)^k \eta_{\wzero}
  \in \re_\wzero, \qquad k \ge 1, \eta \in \re_W.
\ee
If there exists a $k$-th wave, then its result on $\wzero$ is
given by \eqref{e:kthwave}. Applying the above constructions to
$\xi = S^{k-1}_{W}(\eta)$, we obtain that the $k$-th wave
(if there is one) is represented by $\cT_{W} (S^{k-1}_{W}(\eta))$.

If now $\xi \in \re_\wzero$ is distributed according to $\mu_\wzero$,
then $\cT_{W}(\xi)$ is a random subtree of $\Zd$. We will prove
that this random tree has a weak limit $\cT$, which
is almost surely finite.
But first let us show that this is actually
sufficient for finiteness of all waves.

Consider the first wave, and let $W \supseteq V$. By construction,
$\Xi_{W}^1(\eta)$ is precisely the vertex set of
$\cT_{W}(\eta_\wzero)$, therefore
\eqn{e:tree}
{ {\mu_W} ( \Xi_{W}^{1}(\eta) \not\subseteq V )
  = {\mu_W} ( \cT_{W}(\eta_\wzero) \not\subseteq V ). }
Here the right hand side is determined by the distribution
of $\eta_\wzero$ under $\mu_W$. This is different from the law
of $\eta_\wzero$ under $\mu_\wzero$, which is simply the uniform
measure on $\cR_\wzero$. It is the latter that we can get
information about using the correspondence to spanning trees.
Indeed, under $\mu_\wzero$, the spanning tree corresponding to
$\eta_\wzero$ is uniformy distributed on the set of spanning
trees of $\hW$.
In order to translate our results back to $\mu_W$, we show that
the former distribution has a bounded density with respect to
the latter. This will be a consequence of the following lemma.

\begin{lemma}
\label{lem:Re-bnd}
Assume $d \geq 3$. There is a constant $C(d) >0$ such that
\be
\label{e:Re-ratio}
  \sup_{V \subset \Zd} \frac{|\re_{V\setminus\{0\}}|}{|\re_V|}
  \leq C(d),
\ee
where the supremum is over finite sets.
\end{lemma}
\bpr
By Dhar's formula \eqref{Dharformula1},
\[
  |\re_{V\setminus\{0\}}|
  = \det (\Delta_{V\setminus\{0\}})
  = \det (\Delta'_V)
\]
where $\Delta'_V$ denotes the matrix indexed by sites
$y \in V$ and defined by
$(\Delta'_V)_{yz} = (\Delta_{V\setminus\{0\}})_{yz}$
for $y,z \in V \setminus \{ 0 \}$, and
$(\Delta'_V)_{0z} = (\Delta'_V)_{z0} = \delta_{0}(z)$.
We have
\[
  \Delta_{V} + P
  = \Delta'_V
\]
where $P$ is a matrix for which $P_{yz} = 0$
unless $y,z \in N = \{ u \in \Zd : |u| \le 1 \}$.
Moreover, $\max_{y,z \in V} |P_{yz}| \leq 2d-1$. Hence
\[
  \frac{|\re_{V \setminus \{0\}}|}{|\re_V|}
  = \frac{\det(\Delta_V +P)}{\det (\Delta_V)}
  = \det ( I + G_V P).
\]
We have $(G_V P)_{yz} = 0$ unless
$z \in N$. Therefore
\be
\label{kast}
  \det(I + G_V P) = \det \left( (I + G_V P)_N \right).
\ee
By transience of the simple random walk in $d \geq 3$, we have
$\sup_{V} \sup_{y,z} G_V (y,z) \leq G(0,0) < \infty$, and therefore
the determinant of the matrix
$(I + G_V P)_N$ in (\ref{kast}) is bounded by
a constant depending on $d$.
\epr

We note that an alternative proof of Lemma \ref{lem:Re-bnd} can
be given based on the following idea. Consider the graph $\bar{W}$
obtained by adding an extra edge $e$ between $0$ and $\delta_W$ in
$\tW$. Then the ratio in \eqref{e:Re-ratio} can be expressed in
terms of the probability that a uniformly chosen spanning tree
of $\bar{W}$ contains $e$. By standard spanning tree results
\cite[Theorem 4.1]{blps01}, the latter is the same as the probability
that a random walk in $\bar{W}$ started at $0$ first hits $\delta_W$
through $e$.

We continue with the bounded density argument.
For any configuration $\xi \in \re_\wzero$ we have
\be
  {\mu_W} (\eta_\wzero = \xi)
  = \frac{1}{|\re_W|}
    \left| \{ k \in\{ 1,\ldots,2d\} : (k)_0 \xi_\wzero \in \re_W \right|.
\ee
Therefore,
\be
  \frac{{\mu_W} (\eta_\wzero = \xi)}
    {{\mu_\wzero} (\eta_\wzero = \xi)}
  \leq \frac{|\re_\wzero|}{|\re_W|} 2d
  \leq C,
\ee
where, by \eqref{e:Re-ratio}, $C > 0$ does not depend on $\xi$
or on $W$. From this estimate, it follows that
\eqn{e:bnd-ratio}
{  \frac{{\mu_W} (\cT_{W}(\eta_\wzero) \not\subset V)}
     {{\mu_\wzero} (\cT_{W}(\eta_\wzero) \not\subset V)}
   \leq C. }

For a more convenient notation, we let $\nu^{(0)}_W$ denote the
probability measure assigning equal mass to each spanning tree of
$\hW$, or alternatively, to each two-component spanning tree of
$\tW$. We can view $\nu^{(0)}_W$ as a measure on $\{0,1\}^\Ed$ in a
natural way, where $\Ed$ is the set of edges of $\Zd$.
By the Majumdar-Dhar correspondence \cite{MD92}, $\nu^{(0)}_W$
corresponds with the measure $\mu_\wzero$, and the law of
$\cT_{W}$ under $\mu_\wzero$ is that of the component of $0$
under $\nu^{(0)}_W$. We keep the notation $\cT_{W}$ when
referring to $\nu^{(0)}_W$.

We are ready to present the proof of Proposition \ref{prop:finite}
based on the proposition below, whose proof is given in
Section \ref{sec:finiteness}.

\bp
\label{prop:weaklimit}
\begin{itemize}
\item[(i)] For any $d \ge 1$, the limit
  $\lim_W \nu^{(0)}_W = \nu^{(0)}$ exists.
\item[(ii)] Assume $d \geq 3$. Denote the component of $0$
  under $\nu^{(0)}$ by $\cT$. Then
  $\nu^{(0)} (|\cT| < \infty) = 1$.
\end{itemize}
\ep

\emph{Proof of Proposition \ref{prop:finite}
(assuming Proposition \ref{prop:weaklimit}).}
By Proposition \ref{prop:weaklimit} (i), we have
\eqnspl{e:weaklimit}
{ \lim_{W \supseteq V} {\mu_\wzero}
  ( \cT_{W}(\eta_\wzero) \not\subseteq V )
  = \lim_{W \supseteq V} \nu^{(0)}_W ( \cT_W \not\subseteq V)
  = \nu^{(0)} ( \cT \not\subseteq V ). }
By Proposition \ref{prop:weaklimit} (ii), the right hand side of
\eqref{e:weaklimit} goes to zero as $V \uparrow \Zd$, and together
with \eqref{e:bnd-ratio} and \eqref{e:tree}, we obtain the $i = 1$
case of \eqref{lap}.

Finiteness of the other waves follows similarly. For $k \ge 2$
we have by \eqref{e:bnd-ratio}
\eqnspl{e:kth}
{ {\mu_W} ( \Xi_{W}^k (\eta) \not\subseteq V )
  &\le {\mu_W} ( \cT_{W}(S^{k-1}_W \eta) \not\subseteq V ) \\
  &\le C {\mu_\wzero} (\cT_{W}(S^{k-1}_W \eta) \not\subset V) \\
  &= C {\mu_\wzero} ( \cT_{W} (\eta_\wzero) \not\subseteq V ), }
where the last step follows by invariance of $\mu_\wzero$
under $\prod_{j \sim 0} a_{j,\wzero}$. We have already seen
in \eqref{e:weaklimit}
that the right hand side of \eqref{e:kth} goes to zero
as $W, V \uparrow \Zd$, which completes the proof of
Proposition \ref{prop:finite}.

\section{Finiteness of two-component spanning trees}
\label{sec:finiteness}

In this section, we complete the arguments for finiteness of avalanches
by proving Proposition \ref{prop:weaklimit}, which amounts to showing
that the weak limit of $\cT_W$ as $W \uparrow \Zd$ is almost surely
finite. For this, we briefly review below some results on the uniform
spanning forest; see \cite{blps01,lyons05} for more background.

The main statement of Proposition \ref{prop:weaklimit} is part (ii).
It can be deduced from a well-known theorem, namely that all
trees in the uniform spanning forest in $\Zd$ have a single end.
This fact has been known for more general graphs than $\Zd$,
see \cite[Theorem 10.1]{blps01}.
Russell Lyons informed us (private communication) that he,
Ben Morris and Oded Schramm have independently proved a more
general version of statement (ii) in the context of giving a new
and more widely applicable proof of the single end theorem;
see \cite{lms05} and \cite[Chapter 9]{lyons05}.

For finite $W \subset \Zd$, let $\nu_W$ denote the probability
measure assigning equal weight to each spanning tree of $\tW$.
$\nu_W$ is known as the uniform spanning tree measure in $W$
with wired boundary conditions (UST).
We use the algorithm below, due to Wilson \cite{Wilson}, to analyze
random samples from $\nu^{(0)}_W$ and $\nu_W$.

Let $G$ be a finite connected graph. By simple random walk on $G$ we
mean the random walk which at each step jumps to a random neighbour,
chosen uniformly. For a path $\pi = [\pi_0, \dots, \pi_m]$
on $G$, define the loop-erasure of $\pi$, denoted $\LE(\pi)$,
as the path obtained by erasing loops chronologically from $\pi$.

\emph{Wilson's algorithm.}
Pick a vertex $r \in G$, called the root. Enumerate the vertices
of $G$ as $x_1, \dots, x_k$. Let  $(S^{(i)}_n)_{n \ge 1}$,
$1 \le i \le k$ be independent simple random walks started at
$x_1, \dots, x_k$, respectively. Let
\eqnst
{ T^{(1)} = \min \{ n \ge 0 : S^{(1)}_n = r \}, }
and set
\eqnst
{ \gamma^{(1)} = \LE( S^{(1)}[0, T^{(1)}] ). }
Now recursively define $T^{(i)}$, $\gamma^{(i)}$, $i = 2, \dots, k$
as follows. Let
\eqnst
{ T^{(i)}
  = \min \{ n \ge 0 : S^{(i)}_n \in \cup_{1 \le j < i} \gamma^{(j)} \}, }
and
\eqnst
{ \gamma^{(i)} = \LE( S^{(i)}[0, T^{(i)}] ). }
(If $x_i \in \cup_{1 \le j < i} \gamma^{(j)}$, then $\gamma^{(i)}$ is
the single point $x_i$.) Let $T = \cup_{1 \le i \le k} \gamma^{(i)}$.
Then $T$ is a spanning tree of $G$ and is uniformly distributed
\cite{Wilson}.

Applying the algorithm with $G = \tW$ and root $\delta_W$ gives
a sample from $\nu_W$. Similarly, applying the method
with $G = \hW$ and root $\delta_\hW$ we get a sample from $\nu^{(0)}_W$.
It will be convenient to think of the latter construction also taking
place in $\tW$, via the one-to-one correspondence between the
edges of $\hW$ and $\tW$. Note that under this correspondence, a path in
$\hW$ that does not use $\delta_\hW$ as an internal vertex, maps
to a path in $\tW$. Hence the two-component spanning tree in $\tW$
can be built from loop-erased random walks by regarding
$\{ 0, \delta_W \}$ as the ``root''. In other words, the
walks attach either to a piece growing from $0$, or to a piece growing
from $\delta_W$, and these two growing pieces yield the two components.

One can extend Wilson's algorithm to infinite graphs $G$ if random walk
on $G$ is transient \cite{blps01}.
In this case, one chooses the root to be "at infinity", and note that
loop-erasure makes sense for infinite paths that visit each site
finitely many times.

The measures $\nu_W$ can be realized on the same sample space,
$\{0,1\}^\Ed$, as $\nu^{(0)}_W$ introduced
earlier. It is well known that $\nu_W$ has a weak limit
$\nu$ as $W \uparrow \Zd$, called the (wired) uniform spanning
forest (USF) on $\Zd$ \cite{pmnt91,blps01}.
When $d \ge 3$, the USF can be constructed
directly by Wilson's method in $\Zd$, rooted at infinity
\cite[Theorem 5.1]{blps01}.

We write $\omega$ for the random set of edges present under
$\nu$, that is, we identify $\omega \in \{0,1\}^\Ed$ with
the set of edges $e$ for which $\omega(e) = 1$. This allows
us to view $\omega$ as a (random) subgraph of $\Ed$.
We say that an infinite tree $T$ has one end, if there are
no two disjoint infinite paths in $T$. It is known that
\eqnsplst
{ \nu( \text{all components of $\omega$ are infinite trees} )
  &= 1, \\
  \nu( \text{each component of $\omega$ has one end} )
  &= 1; }
see \cite{blps01,lms05}.

\emph{Proof of Proposition \ref{prop:weaklimit}.}
Denote the random set of edges present in the two-component
spanning tree of $\tW$ by $\omega_W$.
Let $W_n$ be an increasing sequence of finite volumes
exhausting $\Zd$. If $B$ is a finite set of edges,
\cite[Corollary 4.3]{blps01} implies that
$\nu^{(0)}_{W_n}(B \subset \omega_{W_n})$
is increasing in $n$. This is sufficient to imply the
weak convergence $\lim_W \nu^{(0)}_W = \nu^{(0)}$,
and the limiting spanning forest
$\omega$ is uniquely determined by the conditions
\eqnst
{ \nu^{(0)} (B \subset \omega)
  = \lim_{n \to \infty} {\nu^{(0)}_{W_n}} (B \subset \omega_{W_n}), }
as $B$ varies over finite edge-sets (see the discussion in
\cite[Section 5]{blps01}).
This proves part (i) of the proposition.

For part (ii), assume $d \ge 3$. The configuration under
$\nu^{(0)}$ can be constructed by Wilson's method directly, by
\cite[Theorem 5.1]{blps01}. Since here $0$ is part of the boundary,
the simple random walks in this construction are either killed when they
hit the component growing from $0$, or they attach to a component
growing from infinity.

Assume now that $\nu^{(0)} ( |\cT| = \infty ) = c_1 > 0$,
and we reach a contradiction.
We consider the construction of the configuration under $\nu^{(0)}$
via Wilson's algorithm. Suppose that the first random walk, call it
$S^{(1)}$, starts from $x \not= 0$.
Write $x \cnctd y$ to denote that $x$ and $y$
are in the same component. Then we have
\eqnst
{ \nu^{(0)} ( x \ntcnctd 0 )
  = \Prob ( \text{$S^{(1)}$ does not hit $0$} )
  = 1 - \frac{G(x,0)}{G(0,0)}
  \to 1 \qquad \text{as $|x| \to \infty$.} }
In particular, there exists an $x \in \Zd$, such that
\eqn{e:positive}
{ \nu^{(0)} ( |\cT| = \infty ,\, x \ntcnctd 0 )
  \ge c_1 / 2. }
Fix such an $x$. Let $B(x,n)$ denote the box of radius $n$ centered at
$x$. Fix $n_0$ such that $0 \in B(x,n_0)$, and $0$ is not a boundary
point of $B(x, n_0)$. By inclusion of events, \eqref{e:positive} implies
\eqn{e:positive2}
{ \nu^{(0)} ( 0 \cnctd \partial B(x,n) ,\, x \ntcnctd 0 )
  \ge c_1 / 2 }
for all $n \ge n_0$. For fixed $n \ge n_0$, let
$y_1 = x$, and let $y_2, \dots, y_K$ be an enumeration of the sites
of $\partial B(x,n)$. We use Wilson algorithm with this enumeration of
sites. Let $S^{(i)}$ and $T^{(i)}$ denote the $i$-th random walk and
the corresponding hitting time determined by the algorithm. We use these
random walks to analyze the configuration under both $\nu^{(0)}$ and
$\nu$.

The event on the left hand side of \eqref{e:positive2} can be recast as
\eqn{e:stops}
{ \{ T^{(1)} = \infty ,\, \text{$\exists\, 2 \le j \le K$ such that
     $T^{(j)} < \infty ,\, S^{(j)}_{T^{(j)}} = 0$} \}, }
and hence this event has probability at least $c_1 / 2$. On the above
event, there is a first index $N$, $2 \le N \le K$, such that the walk
$S^{(N)}$ hits $B(x,n_0)$ at some random time $\sigma$, where
$\sigma < T^{(N)}$. Let $A$ be the latter event. Since $A$ contains
the event in \eqref{e:stops}, $A$ also has probability at least
$c_1 / 2$. Let $p = p (x, n_0)$ denote the
minimum over $z \in \partial B(x, n_0)$ of the probability that a
random walk started at $z$ hits $x$ before $0$ without exiting
$B(x, n_0)$. Clearly, $p > 0$, and is independent of $n$.

Let $B$ denote the subevent of $A$ on which after time $\sigma$, the
walk $S^{(N)}$ hits the loop-erasure of $S^{(1)}$ before hitting $0$
(and without exiting $B(x, n_0)$). We have $\Prob ( B \,|\, A ) \ge p$.
Now we regard the random walks as generating $\nu$. By the definition
of $N$, on the event $A \cap B$, the hitting times
$T^{(1)}, \dots, T^{(N)}$, have the same values as in the construction
for $\nu^{(0)}$, since the walks do not hit $0$.
Moreover, on $A \cap B$, the tree containing $x$ has two disjoint paths
from $\partial B(x, n_0)$ to $\partial B(x, n)$: one is part of the
infinite path generated by $S^{(1)}$, the other part of the path
generated by $S^{(N)}$. Therefore, the probability of the existence
of two such paths is at least $p (c_1 / 2)$, for all $n \ge n_0$.
However, this probability should go to zero as $n \to \infty$,
because under $\nu$, each tree has one end almost surely. This is a
contradicion, proving part (ii) of the proposition.

\section{Tail triviality of $\mu$}
\label{sec:tailtriv}
In this section we study ergodic properties of $\mu$.
In Section \ref{sec:ergodicity} we are going to use the
$d \geq 3$ part of the following theorem.

\bt
\label{thm:tail}
The measure $\mu$ is tail trivial for any $d \ge 2$.
\et

Our proof of Theorem \ref{thm:tail} is divided into two parts.
The argument in the case $2 \le d \le 4$ is quite simple, and is
given in Section \ref{ssec:2d4}. The case $d > 4$ is quite involved,
and is given in Sections \ref{ssec:coding} and \ref{ssec:d4}.

\subsection{The case $2 \le d \le 4$}
\label{ssec:2d4}

\emph{Proof of Theorem \ref{thm:tail} [Case $2 \le d \le 4$].}
The proof is based on the fact that the uniform spanning
forest measure $\nu$ is tail trivial \cite[Theorem 8.3]{blps01}.
Let $\cX \subset \{0,1\}^\Ed$ denote the set of spanning trees of
$\Zd$ with one end. Recall the uniform spanning forest measure
$\nu$ from Section \ref{sec:finiteness}. It was shown by
Pemantle \cite{pmnt91} that when $2 \le d \le 4$, the measure
$\nu$ is concentrated on $\cX$. We can regard any
$\omega \in \cX$ as a tree "rooted at infinity", that is, we
call $x$ an ancestor of $y$ if and only if $x$ lies on the
unique path from $y$ to infinity.

It is shown in \cite{aj} that there is a mapping
$\psi : \cX \to \Omega$ such that $\mu$ is the image of $\nu$
under $\psi$. Moreover, $\psi$ has the following property.
Let $T_x = T_x(\omega)$ denote the tree consisting of all ancestors
of $x$ and its $2d$ neighbours in $\omega$. In other words, $T_x$ is
the union of the paths leading from $x$ and its neighbours to infinity.
It follows from the results in \cite{aj} that
$\eta_x = (\psi(\omega))_x$ is a function of $T_x$ alone.

Assume that $f(\eta)$ is a bounded tail measurable function. Then
for any $n$, $f$ is a function of $\{ \eta_x : \|x\|_\infty \ge n \}$
only. This means that $f(\eta) = f(\psi(\omega)) = g(\omega)$ is a
function of the family $\{ T_x(\omega) : \|x\|_\infty \ge n \}$. Let
$\cF_k = \sigma( \omega_e : e \cap [-k,k]^d = \emptyset )$.
For $1 \le k < n$ consider the event
\eqnst
{ E_{n,k}
  = \bigcap_{x : \|x\|_\infty \ge n}
    \{ T_x \cap [-k,k]^d = \es \}. }
Observe that $E_{n,k} \in \cF_k$, and $g I[E_{n,k}]$ is
$\cF_k$-measurable. Using that $\omega$ has a single end
$\nu$-a.s., it is not hard to check that for any $k \ge 1$
\eqnst
{ \lim_{n \to \infty} \nu (E_{n,k}) = 1. }
Letting $n \to \infty$, this implies that there is an
$\cF_k$-measurable function $\hg_k$, such that $g = \hg_k$
$\nu$-a.s. Since this holds for any $k \ge 1$, tail triviality
of $\nu$ implies that $g$ is constant $\nu$-a.s.
Letting $c$ denote the constant, this implies
\eqnst
{ \mu( f(\eta) = c ) = \nu( f(\psi(\omega)) = c ) = 1, }
which completes the proof in the case $2 \le d \le 4$.
\qed

\subsection{Coding of the sandpile in the case $d > 4$}
\label{ssec:coding}

The simple proof in Section \ref{ssec:2d4} does not work when
$d > 4$, due to the fact that the coding of the sandpile
configuration by the USF breaks down. Nevertheless, it turns out
that a coding is possible if we add extra randomness to the USF,
namely, a random ordering of its components.
Due to the presence of this random ordering,
however, we have not been able to deduce tail triviality of $\mu$
directly from tail triviality of $\nu$, and we need a separate argument.

We start by recalling results from \cite{aj}. Let $\cX$ denote
the set of spanning forests of $\Zd$ with infinitely many components,
where each component is infinite and has a single end. The USF
measure $\nu$ is concentrated on $\cX$ \cite{blps01}. Given
$x \in \Zd$ and $\omega \in \cX$, let
$T^{(1)}_x(\omega), \dots, T^{(k)}_x(\omega)$ denote the trees
consisting of all ancestors of $x$ and its $2d$ neighbours in
$\omega$. Here $k = k_x(\omega) \ge 1$. Each $T^{(i)}_x$
is a union of infinite paths starting at $x$ or a neighbour of $x$,
and has a unique vertex $v^{(i)}_x$ where these paths "first meet".
In other words, $v^{(i)}_x$ is the first vertex that is common to
all of the paths. Let $F^{(i)}_x(\omega)$ denote the finite tree
consisting of all descendants of $v^{(i)}_x$ in $T^{(i)}_x(\omega)$.
Let $\cF$ denote the collection of finite rooted trees in $\Zd$.
Let $\Sigma_l$ denote the set of permutations of the symbols
$\{1, \dots, l\}$.

The sandpile height at $x$ is a function of
$\{ F^{(i)}_x(\omega), v^{(i)}_x(\omega) \}_{i=1}^k$ and a random
$\sigma_x \in \Sigma_k$, in the following sense.

\begin{lemma}
There are functions $\psi_l : \cF^l \times \Sigma_l$, $l = 1, 2, \dots$
such that if $\sigma_x$ is a uniform random element of
$\Sigma_k$, given $\omega$, then
\eqn{e:one-dim-dist}
{ \eta_x = \psi_{k_x} ( (F^{(1)}_x, v^{(1)}_x), \dots,
  (F^{(k)}_x, v^{(k)}_x), \sigma_x ) }
has the distribution of the height variable at $x$ under $\mu$.
\end{lemma}

\begin{proof}
This follows from the proofs of Lemma 3 and Theorem 1 in \cite{aj}.
\end{proof}

\begin{remark}
Here it is convenient to think of $\sigma_x$ as a random ordering
of those components of $\omega$ that contain at least one neighbour
of $x$. Then one can also view $\eta_x$ as a function of
$\{ T^{(i)}_x \}_{i=1}^k$ and $\sigma_x$.
\end{remark}

Next we turn to a description of the joint distribution of
$\{ \eta_x \}_{x \in A_0}$ for finite $A_0 \subset \Zd$.
Let $A = A_0 \cup \partial_e A_0$. Let $\cC^{(1)}, \dots, \cC^{(K)}$,
$K = K_A(\omega)$, denote the components of the USF intersecting
$A$. Each $\cC^{(i)}$ contains a unique vertex $v^{(i)}_A$ where
the paths from $A \cap \cC^{(i)}$ to infinity first meet.
Let $F^{(i)}_A$ denote the finite tree consisting in the portion of
these paths up to $v^{(i)}_A$. In other words, $F^{(i)}_A$ is the
union of the paths from $A \cap \cC^{(i)}$ to $v^{(i)}_A$.
Each rooted tree $(F^{(j)}_x, v^{(j)}_x)$,
$x \in A_0$, $1 \le j \le k_x$ is a subtree
of some $F^{(i)}_A$, $1 \le i \le K$ and the former are determined
by the latter. Let $\sigma_A \in \Sigma_K$ be uniformly distributed,
given $\omega$. For each $x \in A_0$, $\sigma_A$
induces a permutation in $\Sigma_{k_x}$, by restriction.
Then the lemma below follows from the results in \cite{aj}.

\begin{lemma}
\label{lem:joint}
The height configuration in $A_0$ is a function of
$\{ (F^{(i)}_A, v^{(i)}_A) \}_{i=1}^K$ and $\sigma_K$. Moreover,
the joint distribution of $\{ \sigma_x \}_{x \in A_0}$ is the one
induced by $\sigma_A$.
\end{lemma}
\qed

From the above, we obtain the following description of $\mu$
in terms of the USF and a random ordering of its components.
Let $\omega \in \cX$ be distributed according to $\nu$.
Given $\omega$, we define a random partial ordering $\prec_\omega$
on $\Zd$ in the following way.
Let $\cC^{(1)}, \cC^{(2)}, \dots$ be an enumeration of the
components of $\omega$, and let $U_1, U_2, \dots$ be i.i.d.~random
variables, given $\omega$, having the uniform distribution on $[0,1]$.
Define the random partial order $\prec_\omega$ depending on $\omega$
and $\{U_i\}_{i \ge 1}$ by letting  $x \prec_\omega y$ if and only if
$x \in \cC^{(i)}$, $y \in \cC^{(j)}$ and $U_i < U_j$. Even though
$\prec_\omega$ is defined for sites, it is simply an ordering of the
components of $\omega$. The distribution of $\prec_{\omega}$ is in fact
uniquely characterized by the property that it induces the uniform
permutation on any finite set of components, and one could define
it by this property, without reference to the $U$'s. This in turn
shows that the distribution is independent of the ordering
$\cC^{(1)}, \cC^{(2)}, \dots$ initially chosen.

Let $\cQ = \{0,1\}^{\Zd \times \Zd}$ denote the space of binary
relations (where for $q \in \cQ$ we interpret $q(x,y) = 1$ as
$x \prec y$, and $q(x,y) = 0$ otherwise). We denote the joint law
of $(\omega, \prec)$ on $\cX \times \cQ$ by $\tnu$.
From the couple $(\omega, \prec)$, we can recover the random
permutations $\sigma_x$ as follows. If $v^{(1)}_x, \dots, v^{(k)}_x$
are as defined earlier, then
\eqn{e:defsigmax}
{ (\sigma_x(1), \dots, \sigma_x(k)) = (j_1, \dots, j_k)
    \qquad \text{if and only if} \qquad
    v^{(j_1)}_x \prec_\omega \dots \prec_\omega v^{(j_k)}_x. }
The discussion above, and Lemma \ref{lem:joint} easily implies
the following lemma.

\begin{lemma}
\label{lem:psiexists}
Suppose that $(\omega, \prec_\omega)$ has distribution $\tnu$.
Let $\eta_x$ be given by \eqref{e:one-dim-dist}, where
$\sigma_x$ is defined by \eqref{e:defsigmax}.
Then $\{ \eta_x \}_{x \in \Zd}$ has distribution $\mu$. In
particular, there is a $\tnu$-a.s.~defined function
$\psi: \cX \times \cQ \to \Omega$ such that $\mu$ is the image of
$\tnu$ under $\psi$.
\end{lemma}
\qed

Before we start the argument proper, we need to recall some further
terminology from \cite{aj}. Given finite rooted trees
$(\bF, \bv) = (F_i, v_i)_{i=1}^K$ and a finite set of sites $A$,
define the events
\eqnsplst
{ D(\bv) &= \{ \text{$v_1, \dots, v_K$ are in distinct components
    of $\omega$} \}, \\
  B(\bF, \bv) &= D(\bv)
    \cap \{ \text{$F^{(i)}_A = F_i, v^{(i)}_A = v_i$ for
    $1 \le i \le K$} \}, }
We also need versions of these events for finite $\La \subset \Zd$.
The wired UST $\omega_\La$ in volume $\La$ can be viewed as the
union of one or more components ($x, y \in \La$ are in the same
component if they are connected without using the special vertex
$\delta_\La$). Let
$\cC^{(1)}_\La, \dots, \cC^{(K)}_\La$ be the list of components
intersecting $A$. We define $v^{(i)}_{A,\La}$
and $F^{(i)}_{A,\La}$ analogously to the infinite volume case,
this time using the components $\cC^{(i)}_\La$. Now we define
\eqnsplst
{ D_\La(\bv) &= \{ \text{$v_1, \dots, v_K$ are in distinct components
    of $\omega_\La$} \}, \\
  B_\La(\bF, \bv) &= D_\La(\bv)
    \cap \{ \text{$F^{(i)}_{A,\La} = F_i, v^{(i)}_{A,\La} = v_i$ for
    $1 \le i \le K$} \}. }

\subsection{Proof in the case $d > 4$}
\label{ssec:d4}

\subsubsection{Outline of the proof}
\label{sssec:outline}

Recall that tail triviality is equivalent to the following
\cite[Proposition 7.9]{georgii}. For any cylinder event $E'$ and $\eps > 0$
there exists $n$ such that (with $V_n = [-n,n]^d \cap \Zd$) for any
event $R' \in \cF_{V_n^c}$ we have
\eqn{e:equiv-tail}
{ | \mu(E' \cap R') - \mu(E') \mu(R') | \le \eps. }
Let $E = \psi^{-1}(E')$ and $R = \psi^{-1}(R')$, where $\psi$ is as in
Lemma \ref{lem:psiexists}. Suppose that $E'$
depends on the sites in the finite set $A_0$, and put
$A = A_0 \cup \partial_e A_0$.

For the proof, we want to show that $E$ and $R$ ``decouple'', if
$n$ is sufficiently large. We try to achieve this by showing that
they can be approximated by events that depend on portions of $\omega$
that are spatially separated. The main difficulty is that dependence
between $E$ and $R$ also exists due to the ordering $\prec$, and it
requires work to show that the dependence on $\prec$ also decouples.
Below we give a rough outline of strategy for this.

By Lemma \ref{lem:joint}, the occurrence or not of the event $E$ is
determined by a collection of finite tree subgraphs $(F^{(i)}_A, v^{(i)}_A)$
of $\omega$, and an ordering of these trees. We get an approximation of the
event $E$, if we consider the contribution of only those configurations
for which $F^{(i)}_A \subset V_r$, $1 \le i \le K$ for some large $r$.
Suppose that we have also approximated $R$ by an event that depends
on the restriction of $\omega$ to $V_m^c$, where $r < m < n$.
Condition on the restriction of $\omega$ to $V_m^c$, and also on the
restriction of $\prec$ to this portion of $\omega$. In particular,
for any $w_1, w_2 \in \partial V_m$, our conditioning specifies whether
$w_1 \prec w_2$ or not.

For fixed $1 \le i < j \le K$, let $\pi_i$ and $\pi_j$ denote
the paths in $\omega$ from $v^{(i)}_A$ and $v^{(j)}_A$, respectively,
to $\partial V_m$. Suppose $\pi_i$ and $\pi_j$ end in vertices $w(i)$ and $w(j)$,
respectively. The conditional probability of $\{ v^{(i)}_A \prec v^{(j)}_A \}$
is determined by the conditional probability of $w(i) \prec w(j)$.
If $r \ll m$, then due to fluctuations in the behaviour of the paths
$\pi_i$ and $\pi_j$, the conditional probability of the
events $\{ w(i) \prec w(j) \}$ and $\{ w(j) \prec w(i) \}$
will be approximately equal, and we obtain the desired decoupling.

In the next section, we specify suitable approximations of
$E$ and $R$.

\subsubsection{Approximating $E$ and $R$}
\label{sssec:approx}

We first have a closer look at the event $E$. We define
\eqnsplst
{ S(\bF, \bv, \sigma) &= B(\bF, \bv)
   \cap \{ v_{\sigma(1)} \prec \dots \prec v_{\sigma(K)} \}, \\
  \cG_E &= \{ (\bF, \bv, \sigma) : S(\bF, \bv, \sigma) \subset E \}, \\
  \cG_E(r) &= \{ (\bF, \bv, \sigma) \in \cG_E :
    \text{$F_i \subset V_r$ for $1 \le i \le K$} \}. }
The event $E$ is a disjoint union of $S(\bF, \bv, \sigma)$ over
$(\bF, \bv, \sigma) \in \cG_E$. By Lemma \ref{lem:psiexists}, we have
\eqn{e:nuE}
{ \mu(E')
  = \tnu(E)
  = \sum_{(\bF, \bv, \sigma) \in \cG_E} \frac1{K!} \nu( B(\bF, \bv) ). }
We also define an analogue of $S$
in a finite volume $\La$. Assume that the relation $\prec_\partial$ is
prescribed on the exterior boundary of $\La$. For any realization
of the wired UST $\omega_\La$ there is a unique extension of
$\prec_\partial$ into $\La$, denoted $\prec_\La$, where $x \prec_\La y$
if and only if they are connected (in $\omega_\La$) to boundary
vertices $w(x)$ and $w(y)$ satisfying $w(x) \prec_\partial w(y)$.
Using this extension, we define
\eqnst
{ S_\La(\bF, \bv, \sigma) = B_\La(\bF, \bv)
    \cap \{ v_{\sigma(1)} \prec_\La \dots \prec_\La v_{\sigma(K)} \}. }
We let $\tnu_{\La, \prec_\partial}$ denote the law of
$(\omega_\La, \prec_\La)$ with boundary condition $\prec_\partial$.

Introduce
\eqnst
{ G = G(r) = \{ \text{$F^{(i)}_A \subset V_r$ for $1 \le i \le K$} \}, }
where we asume that $A_0 \subset V_r \subset V_n$. Now $E \cap G$
is a disjoint union of the events $S(\bF, \bv, \sigma)$ over
$(\bF, \bv, \sigma) \in \cG_E(r)$. Since $A$ is fixed, there exists
$r_0(\eps)$ such that for $r \ge r_0(\eps)$ we have $\nu(G(r)^c) \le \eps$.
The event $E \cap G(r)$ will serve as an approximation for $E$.

Turning to $R$, we define
\eqnsplst
{ \cH &= \cH_n = \bigcup_{x \in V_n^c} \bigcup_{i=1}^{k_x}
    \text{ vertex set of $T^{(i)}_x$} \\
  \cD &= \cD_n = \Zd \setminus \cH_n. }
The occurrence of $R$ is determined by the collection of edges joining
vertices in $\cH$ together with the restriction of $\prec$ to $\cH$.
We also introduce for $r < m < n$ and
$V_m \subset \La \subset V_n$ the events
\eqnst
{ F = F(n,m) = \{ \cH_n \cap V_m = \emptyset \}
  \qquad \text{ and } \qquad
  F_\La = F_\La(n,m) = \{ \cD_n = \La \}. }
Here $F$ is the event that the portion of $\omega$
determining the sandpile configuration in $V_n^c$ does not intersect
$V_m$. The event $R \cap F$ will serve as an approximation of $R$,
as mentioned in Section \ref{sssec:outline}. However, we will
further decompose $F$ as the disjoint union
$F = \bigcup_\La F_\La$. The reason is that the conditional law of
$\nu$ inside $\cD$, given $F_\La$ is simple: it is $\nu_\La$
(see Section \ref{sssec:decomp} below).

The value of $m$ will be chosen large with respect to $r$.
It is easy to see that there exists $n_0(m,\eps)$, such that if
$n \ge n_0(m,\eps)$ then $\nu(F^c) \le \eps$. This is because $F(n,m)$
is monotone increasing in $n$, and
$\cap_{n = m+1}^\infty F(n,m)^c = \emptyset$,
since each component of the USF has a single end.

For technical reasons, we will in fact need a further subevent
of $F$, on which, given the configuration in $\cH$,
with high conditional probability:
\eqnst
{ \text{for all $x, y \in V_r$} \quad  x \cnctd y
  \quad \text{implies} \quad
  \text{$x \cnctd y$ inside $V_m$.} }
Let
\eqnst
{ J = \{ \text{$\forall x, y \in V_r$ : if $x \conn y$ then
    $x \conn y$ inside $V_m$} \}. }
There exists $m_0(r, \eps)$, such that if $m \ge m_0(r,\eps)$, then
$\nu(J^c) \le \eps \eps_1$, where we have set
$\eps_1 = \eps_1(r) = \eps / |\cG_E(r)|$. Define the event
\eqnst
{ J_0
  = F \cap \big\{ \nu \big( J^c \,\big|\, \omega_{\cH_n} \big)
    \le \eps_1 \big\}, }
where $\omega_{\cH_n}$ denotes the configuration on the set of edges
touching $\cH_n$. By Markov's inequality,
\eqnst
{ \nu(J_0^c)
  \le \nu(F^c)
      + \nu \big( F \cap \big\{ \nu\big( J^c \,\big|\,
      \omega_{\cH_n} \big) \ge \eps_1 \big\} \big)
  \le \eps + \frac{\nu(J^c)}{\eps_1}
  \le 2 \eps. }

Summarizing the above, if $r \ge r_0$, $m \ge m_0(r,\eps)$ and
$n \ge n_0(m,\eps)$, we have
\eqn{e:likely}
{ | \mu(E' \cap R') - \tnu( E \cap G \cap R \cap J_0 ) ) |
  \le 3 \eps. }

In the next section, we obtain a decomposition of the event
$E \cap G \cap R \cap J_0$, that allows us to analyze it
via Wilson's method.

\subsubsection{Decomposition of $E \cap G \cap R \cap J_0$}
\label{sssec:decomp}

We are going to regard the edges of $\omega$ being directed towards
infinity. By the definition of $\cH$, there are no directed edges
from $\cH$ to $\cD$. Therefore, given the restriction of $\omega$
to $\cH$, the
conditional law of $\omega$ in $\cD$ is that of the wired uniform
spanning tree in $\cD$, that is $\nu_D$. One can see this by
using Wilson's method rooted at infinity to first generate
$\cH$ and the configuration on $\cH$, and then the configuration
in $\cD$.

Note that the event $F_\La$ only depends on the portion of $\omega$
outside $\La$. We want to rewrite the second term on the left hand side
of \eqref{e:likely} by conditioning on $F_\La$, the portion of $\omega$
outside $\La$, and the restriction of $\prec$ to $\Zd \setminus \La$.
By the previous paragraph, the conditional distribution of
$(\omega, \prec)$ inside $\La$ is given by $\tnu_{\La, \prec_\partial}$,
where $\prec_\partial$ is determined by the conditioning.

The above implies
\eqn{e:D-decomp}
{ \tnu( E \cap G \cap R \cap J_0 )
  = \sum_{V_m \subset \La \subset V_n}
    \int_{R \cap J_0 \cap F_\La}
    \tnu_{\La, \prec_\partial} ( E \cap G )\, d\tnu. }
Since the integration in \eqref{e:D-decomp} is over a subset of $J_0$,
in what follows, we assume that the boundary condition
$\prec_\partial$ is compatible with the event $J_0$, in the sense that
it arises from a configuration belonging to $J_0$.
The expression $\tnu_{\La, \prec_\partial} ( E \cap G )$
can be further decomposed as follows:
\eqn{e:elem}
{ \tnu_{\La, \prec_\partial} ( E \cap G )
  = \sum_{(\bF, \bv, \sigma) \in \cG_E(r)}
    \tnu_{\La, \prec_\partial} ( S_\La (\bF, \bv, \sigma) ). }

In the remainder of the proof our aim is to show that the summand
in \eqref{e:elem} is close to $\nu( B(\bF, \bv) ) / K!$,
uniformly in $\Lambda$ and the boundary condition,
if $m$ is large enough. In the next section, we formulate
precisely the statement we need as Lemma \ref{lem:decoup},
and prove the theorem given Lemma \ref{lem:decoup}.
Finally, in Section \ref{sssec:prooflemma}, we complete
the argument by proving Lemma \ref{lem:decoup}.

\subsubsection{Decoupling lemma and proof of theorem}
\label{sssec:prooftheorem}

\begin{lemma}
\label{lem:decoup}
There exists a universal constant $C$ and $m_1(r,\eps)$, such that
for any $(\bF, \bv, \sigma) \in \cG_E(r)$,
$m \ge \max\{ m_0(r,\eps), m_1(r,\eps)\}$,
$V_m \subset \La$ and any boundary condition
$\prec_\partial$ compatible with $J_0$, we have
\eqn{e:decoup}
{ \left| \tnu_{\La, \prec_\partial} ( S_\La (\bF, \bv, \sigma) )
    - \frac{\nu( B(\bF, \bv) )}{K!} \right|
  \le C \eps_1 = C \frac{\eps}{|\cG_E(r)|}. }
\end{lemma}

\emph{Proof of Theorem \ref{thm:tail} [Case $d > 4$] assuming
Lemma \ref{lem:decoup}.}
Given $\eps > 0$ let $r \ge r_0(\eps)$,
$m \ge \max\{ m_0(r,\eps), m_1(r,\eps) \}$ and $n \ge n_0(m,\eps)$.
Then the estimate in \eqref{e:decoup}, formula \eqref{e:nuE} and
\eqref{e:elem} imply
\eqnst
{ | \tnu_{\La, \prec_\partial} (E \cap G) - \tnu(E \cap G) |
  \le C \eps. }
Substituting this into \eqref{e:D-decomp}, and performing the
integral and the sum, we get
\eqnst
{ | \tnu (E \cap G \cap R \cap J_0)
    - \tnu(E \cap G) \tnu(R \cap J_0) |
  \le C \eps. }
Due to $r \ge r_0(\eps)$ and $n \ge n_0(m,\eps)$, we have
$\tnu(G^c) \le \eps$ and $\tnu(J_0) \le 2 \eps$, which yields
\eqnst
{ | \mu( E' \cap R') - \mu( E' ) \mu( R' ) | \le C' \eps, }
with a universal constant $C'$. This proves Theorem \ref{thm:tail}
in the case $d > 4$.
\qed

\subsubsection{Proof of decoupling lemma}
\label{sssec:prooflemma}

We prove Lemma \ref{lem:decoup} by analyzing the event $B_\La(\bF, \bv)$
in terms of Wilson's algorithm. The proof is similar to the proof
of Lemma 3 in \cite{aj}, however it does not seem possible to use
that result directly. Before starting the proof proper, we
introduce some notation.

Fix $(F_i, v_i)_{i=1}^K$ and $\sigma \in \Sigma_K$. Let
$A = \{ y_1, \dots, y_{|A|} \}$. We apply
Wilson's method to generate part of the wired UST in $\La$ with the
following enumeration of sites:
\eqnst
{ v_1, \dots, v_K, y_1, \dots, y_{|A|}. }
Let $S^{(i)}$, $i=1,\dots,K$ be independent
simple random walks started at $v_i$. Let $\gamma^{(i)}_\La$ denote
the loop-erasure of $S^{(i)}$ up to its first exit time from $\La$.
We define a random walk event $C_\La$ whose occurrence will be
equivalent to the occurrence of $B_\La(\bF, \bv)$, by Wilson's
method. Since the event $D_\La(\bv)$ has to occur, we require that
for $i = 1, \dots, K$, $S^{(i)}$ upto its first exit time be disjoint
from $\cup_{1 \le j < i} \gamma^{(i)}_\La$. In addition, the fact
that $B_\La (\bF, \bv)$ has to occur, gives conditions on the paths
starting at $y_1, \dots, y_{|A|}$, namely, these paths have to realize
the events $(F^{(i)}_A, v^{(i)}_A) = (F_i, v_i)$, given the paths
$\{ \gamma^{(i)}_\La \}_{i=1}^K$. These implicit conditions define
$C_\La$. More precisely, the loop-erased walk $\eta_1$ started at
$y_1$ has to equal the path in $\cup_i F_i$ from $y_1$ to
$\{ v_1, \dots, v_K \}$. The loop-erased walk $\eta_2$ started at
$y_2$ has to equal the path in $\cup_i F_i$ from $y_2$ to
$\{ v_1, \dots, v_K \} \cup \eta_1$, and so on.

We write $\Prob$ for probabilities associated with random
walk events, and we couple the constructions in different volumes
by using the same infinite random walks $S^{(i)}$.
We also define the random walk event $C$, corresponding to
$B(\bF, \bv)$, analogously to the finite volume case.

\medbreak

\emph{Proof of Lemma \ref{lem:decoup}.}
Let $W^{(i)}_\La$ denote the first vertex $S^{(i)}$ visits in
$\Zd \setminus \La$.
Then we have
\eqn{e:rws}
{ \tnu_{\La, \prec_\partial} ( S_\La (\bF, \bv, \sigma) )
  = \Prob ( C_\La,\, W^{(\sigma(1))}_\La \prec \dots \prec
    W^{(\sigma(K))}_\La ). }
For $r < l < m$ we consider the event $C_{V_l}$, and write
$C_l$ for short. It is not hard to see that
$\lim_\La I[C_\La] = I[C]$, $\Prob$-a.s., which implies that
for $l$ large enough, $\Prob( C_l \triangle C ) \le \eps_1$.
(Here $\triangle$ denotes symmetric difference.)
Hence the difference between the right hand side of \eqref{e:rws} and
\eqn{e:rws-trunk}
{ \Prob ( C_l,\,  W^{(\sigma(1))}_\La \prec \dots \prec
    W^{(\sigma(K))}_\La ) }
is at most $2 \eps_1$. Recall that $\La \supseteq V_m$, and $m > l$.
By conditioning on the first
exit points from $V_l$, \eqref{e:rws-trunk} can be written as
\eqn{e:factored}
{ \Prob( C_l )\, \Prob \big( W^{(\sigma(1))}_\La \prec \dots \prec
    W^{(\sigma(K))}_\La \,\big|\, W^{(1)}_l, \dots, W^{(K)}_l \big). }
The first factor here differs from $\Prob(C) = \nu(B(\bF, \bv))$ by
at most $\eps_1$. If $m$ is large with respect to $l$, the value of
the second factor is essentially independent of $\sigma$. This is
because the distributions of $W^{(i)}_\La$ and $W^{(j)}_\La$ given
$W^{(i)}_l$ and $W^{(j)}_l$ (respectively), can be made arbitrarily
close in total variation distance. This implies that the difference
between \eqref{e:factored} and
\eqnst
{ \Prob(C)\, \Prob \big( W^{(1)}_\La \prec \dots \prec
    W^{(K)}_\La \,\big|\, W^{(1)}_l, \dots, W^{(K)}_l \big) }
is at most $\eps_1$, if $m$ is large enough, uniformly in $\La$.

Observe that if $W^{(i)}_\La \conn W^{(j)}_\La$ for some
$1 \le i < j \le K$, then the event $J^c$ occurs.
Since the boundary condition $\prec_\partial$ is compatible with $J_0$,
we have
\eqn{e:coincide}
{ \Prob \big( C_\La,\, \text{$W^{(i)}_\La \conn W^{(j)}_\La$ for some
     $1 \le i < j \le K$} \big)
  \le \tnu_{\La, \prec_\partial}(J^c)
  \le \eps_1. }
It follows that for some universal constant $C$, if $m$ is large enough
\eqnst
{ \big| \Prob \big( C_\La,\, W^{(\sigma(1))}_\La \prec \dots \prec
    W^{(\sigma(K))}_\La \big) - \Prob(C) / K! \big|
  \le C \eps_1. }
This proves the lemma. \qed

\section{Ergodicity of the stationary process}
\label{sec:ergodicity}
Arrived at this point, we can apply the results
in \cite{mrs}, and  we obtain the following.
\bt
\label{thm:process}
Let $\vi: \Zd \to (0,\infty)$ be an addition rate such that
\be\label{greencon}
\sum_{x \in \Zd} \vi (x) G(0,x)  <\infty.
\ee
Then the following hold.
\ben
\item The closure of the operator on $L_2 (\mu)$ defined
  on local functions
  by
  \be
  L_\vi f = \sum_{x \in \Zd} \vi(x) (a_x-I)f
  \ee
  is the generator of a stationary Markov process $\{ \eta_t: t\geq 0 \}$.
\item   Let $N^\vi_t (x)$ denote Poisson processes with rate $\vi (x)$
  that are independent (for different $x$). The limit
  \be
    \eta_t
    = \lim_{V\uparrow \Zd}
      \left[ \prod_{x \in V} a_x^{N^\vi_t (x)} \right] \eta
  \ee
  exists a.s.~with respect to the product of the
  Poisson process measures on $N^\vi_t$ with
  the stationary measure $\mu$ on the $\eta\in\Omega$.
  Moreover, $\eta_t$ is a cadlag version of the process
  with generator $L_\vi$.
\een
\et

Let $\{ \eta_t: t\geq 0 \}$ be the stationary
process with generator $L_\vi = \sum_x \vi (x) (a_x-I)$.
We recall that a process is called ergodic if
every (time-)shift invariant measurable
set has measure zero or one. For a Markov process,
this is equivalent to the following: if $S_t f = f$
for all $t>0$, then $f$ is constant $\mu$-a.s.
This in turn is equivalent to the statement that $Lf=0$
implies  $f$ is constant $\mu$-a.s.
The tail $\si$-field on $\Omega$ is
defined as usual:
\be
  \fe_\infty = \bigcap_{n\in\N} \si\{ \eta (x): |x|\geq n \}
\ee
A function $f$ is tail measurable if
its value does not change by changing
the configuration in a finite number of
sites, that is, if
\[
  f(\eta) = f(\xi_V \eta_{V^c})
\]
for every $\xi$ and $V\subset\Zd$ finite.
\bt
\label{thm:ergodic}
The stationary process of Theorem \ref{thm:process} is mixing.
\et
\bpr Recall that $G$ denotes the group generated by the unitary
operators $a_x$ on $L_2(\mu)$. Consider the following statements.
\ben
\item The process $\{ \eta_t : t\geq 0 \}$
  is ergodic.
\item The process $\{ \eta_t : t\geq 0 \}$ is
  mixing.
\item Any $G$-invariant function is $\mu$-a.s.~constant.
\item $\mu$ is tail trivial.
\een
Then we have the following implications:
1, 2 and 3 are equivalent and 4 implies 3.
This will complete the proof, because 4 holds by Theorem
\ref{thm:tail}.

It is easy to see that on $L_2 (\mu )$,
\be
  L^* = \sum_{x \in \Zd} \vi (x) (a_x^{-1} - I).
\ee
Hence $L$ and $L^*$ commute, that is, $L$ is a normal
operator. The equivalence of 1 and 2 then follows
from \cite[Lemmas 6 and 7]{rosenblatt} and an
adaptation to continuous time. To see the equivalence of 1 and 3:
invariance of $\mu$ under $a_x$ and $a_x^{-1}$ implies
\beq
  \langle Lf | f \rangle
  &=& -\frac12 \sum_{x \in \Zd} \vi (x)\int (a_x f -f)^2 d\mu,
\nonumber\\
&=&
  \langle L^*f | f \rangle
\nonumber\\
  &=& -\frac12 \sum_{x \in \Zd} \vi (x)\int (a_x^{-1} f -f)^2 d\mu.
\eeq
Hence $Lf = 0$ is equivalent to $f$ being invariant under all
$a_x$ and $a_x^{-1}$, and thus under the action of $G$.
Finally, to prove the implication $4 \Rightarrow 3$, we will
show that a function invariant under the action of $G$ is
tail measurable. Suppose $f:\Omega \to \R$, and
$f = a_x f = a_x^{-1} f$ $\mu$-a.s. for all $x \in \Zd$.
There exists a full measure $G$-invariant subset $\Omega_0$
so that the restriction of $f$ to $\Omega_0$ is $G$-invariant.
If $\eta$ and $\zeta$ are elements of $\Omega_0$ and differ in a
finite number of coordinates, then
\be
  \zeta = \prod_{x \in \Zd} a_x^{\zeta (x) - \eta (x)} \eta
\ee
and hence $f(\eta ) = f(\zeta)$. This implies that $f$ is
$\mu$-a.s.~equal to a tail measurable function.
\epr

\appendix
\section{Appendix}

In this section we show how to extend the argument of
\cite{aj} in the case $d > 4$ and prove $\lim_\La \mu_\La = \mu$.
Using the notation of Section \ref{ssec:coding}, let
\eqnsplst
{ X_{\La, i}
  &= \dist_\La (v_i, \delta_\La), \quad i = 1, \dots, K, \\
  Y_\La
  &= \max_{1 \le i < j \le K} \left| X_{\La, i} - X_{\La, j} \right|. }
where $\dist_\La$ denotes graph distance in the uniform spanning tree
$\omega_\La$. We define the random permutation $\sigma^*_K$ by the
requirement:
\eqnst
{ \sigma^*_K = \sigma \qquad \text{if and only if} \qquad
  X_{\La, \sigma(1)} \le \dots \le X_{\La, \sigma(K)}, }
where we take a fixed but otherwise arbitrary rule to settle ties.
Let $K (\bar F) = \max_{1 \le i \le K} \diam(F_i)$.
The required extension follows once we show the following analogues
of \cite[Eqns. (18) and (19)]{aj}.
\eqn{e:(18)}
{ \lim_{\La} \mu_\La \left(
  B_\La( \bar F, \bar v),\, Y_\La \le K(\bar F) \right)
  = 0, }
and
\eqn{e:(19)}
{ \lim_{\La} \mu_\La \left(
  B_\La( \bar F, \bar v),\, \sigma^*_K = \sigma,\, Y_\La > K(\bar F) \right)
  = \frac{1}{K!} \mu \left( B( \bar F, \bar v) \right). }
Most of the argument in \cite{aj} does apply to general volumes, and here
we detail only those points where differences arise.
We use the notation introduced in Section \ref{sssec:prooflemma}
for Wilson's algorithm.

We start with the proof of \eqref{e:(18)}.
Let $x,y \in \Zd$ be fixed, and let $S^{(1)}$ and
$S^{(2)}$ be independent simple random walks starting at $x$ and
$y$, respectively. Let $T^{(1)}_\La$ and $T^{(2)}_\La$ be the first
exit times from $\La$ for these random walks.
The required extension of \eqref{e:(18)} follows from an extension of
(27)\cite{aj}, which in turn follows from the statement
\eqn{e:noconcent}
{ \lim_{\delta \to 0} \limsup_\La
  \Prob \Big( 1 - \delta \le \frac{T^{(1)}_\La}{T^{(2)}_\La}
    \le 1 + \delta \Big)
  = 0. }
Statement \eqref{e:noconcent} is proved in \cite{jk}.

For the extension of \eqref{e:(19)}, we recall from
Section \ref{ssec:coding} the events $B_\La(\bF,\bv)$ and
$B(\bF,\bv)$ defined for a collection $(F_i, v_i)_{i=1}^K$.
Let $S^{(i)}$, $i=1,\dots,K$ be independent random walks started at
$v_i$, respectively. Let $T^{(i)}_\La$ be the exit time of $S^{(i)}$
from $\La$.
Also recall the random walk events $C_\La$ and $C$, and that
$C_m$ and $T^{(i)}_m$ are short for $C_\La$ and $T^{(i)}_\La$
when $\La = [-m,m]^d \cap \Zd$. By the
arguments in \cite{aj}, the required extension of \eqref{e:(19)}
follows, once we show an extension of (32)\cite{aj}, namely that
for any permutation $\sigma \in \Sigma_K$
\eqn{e:unif1}
{ \lim_{m \to \infty} \lim_\La
  \Prob \left( C_m, T^{\sigma(1)}_\La < \dots < T^{\sigma(K)}_\La \right)
  = \Prob (C) \frac1{K!}. }
Observe that $C_m$ and the collection
$\tT^{(i)}_{\La,m} = T^{(i)}_\La - T^{(i)}_m$, $i = 1, \dots, K$ are
conditionally independent, given $\{ S^{(i)}(T^{(i)}_m) \}_{i=1}^K$.
Therefore, using \eqref{e:noconcent}, the left hand side of
\eqref{e:unif1} equals
\eqn{e:unif2}
{ \lim_{m \to \infty} \lim_\La
  \Prob (C_m)\, \Prob \left( \tT^{\sigma(1)}_{\La,m} < \dots
   < \tT^{\sigma(K)}_{\La,m} \right). }
The second probability approaches $1 / K!$ for any fixed $m$,
and hence the limit in \eqref{e:unif2}
equals $\Prob(C) / K!$. This completes the proof of the required
extension of \eqref{e:(19)}.

\

{\bf Acknowledgements.} We thank Russell Lyons for comments on the
first version of the manuscript, and a referee for comments on the
second version. The work of AAJ was supported by NSERC of Canada, and
partly carried out at the Centrum voor Wiskunde en Informatica, Amsterdam,
The Netherlands.


\begin{thebibliography}{99}

\bibitem{aj} Athreya, S.R.~and J\'arai, A.A.:
  Infinite volume limit for the stationary distribution
  of Abelian sandpile models.
  \CMPsh\ {\bf 249}, 197-213 (2004).
  An erratum for this paper appeared in
  Comm. Math. Phys.  {\bf 264}, 843, (2006)
  with an electronic supplemental material.

\bibitem{btw} Bak, P., Tang, K. and Wiesefeld, K.:
  Self-organized criticality.
  Phys. Rev. A {\bf 38}, 364-374 (1988).

\bibitem{blps01} Benjamini, I., Lyons, R., Peres, Y.~and
  Schramm, O.:
  Uniform spanning forests.
  \AOPsh\ {\bf 29} 1--65 (2001).

\bibitem{Dhar} Dhar D.:
  Self organized critical state of sandpile automaton models.
  Phys.~Rev.~Letters {\bf 64}, 1613-1616  (1990).

\bibitem{Dhar1} Dhar D.:
  The Abelian sandpile and related models.
  Physica A {\bf 263}, 4-25 (1999).

\bibitem{Dhar06} Dhar, D.:
  Theoretical studies of self-organized criticality.
  Physica A, {\bf 369}, 29--70 (2006).

\bibitem{gabr} Gabrielov, A.:
  Asymmetric Abelian avalanches and sandpiles. Preprint 93-65,
  Mathematical Sciences Institute, Cornell University (1993).

\bibitem{georgii} Georgii, H.-O.:
  \emph{Gibbs measures and phase transitions.}
  de Gruyter, Berlin (1988).

\bibitem{IKP94} Ivashkevich, E.V., Ktitarev, D.V.~and Priezzhev, V.B.:
  Critical exponents for boundary avalanches in two-dimensional
  Abelian sandpile.
  J.~Phys.~A {\bf 27} L585--L590 (1994).

\bibitem{IKP94b} Ivashkevich, E.V., Ktitarev, D.V.~and Priezzhev, V.B.:
  Waves of topplings in an Abelian sandpile.
  Physica A {\bf 209} 347--360 (1994).

\bibitem{priz}  Ivaskevich, E.V. and Priezzhev, V.B.:
  Introduction to the sandpile model.
  Physica A {\bf 254}, 97-116 (1998).

\bibitem{Jarai05} J\'arai, A.A.:
  Thermodynamic limit of the Abelian sandpile model on $\Zd$.
  \MPRFsh\ {\bf 11}, 313--336 (2005).

\bibitem{jk} J\'arai, A.A.~and Kesten, H.:
  A bound for the distribution of the hitting time of arbitrary sets
  by random walk.
  \EJPsh\ {\bf 9}, 152--161 (2004).

\bibitem{Lawler} Lawler, G.F.:
  \emph{Intersections of random walks.}
  Birkh\"auser, softcover edition (1996).

\bibitem{Lawler2} Lawler, G.F.:
  Loop-erased random walk.
  In \emph{Perplexing problems in Probability},
  Progress in Probability, Vol.~44, Birkh\"auser, Boston, (1999).

\bibitem{lms05} Lyons, R., Morris, B.~and Schramm, O.:
  Ends in uniform spanning forests. In preparation.

\bibitem{lyons05} Lyons, R. and Peres, Y.:
  \emph{Probability on trees.}
  Book in preparation.
  {\tt http://mypage.iu.edu/\string~rdlyons/prbtree/prbtree.html.}

\bibitem{mrs} Maes, C., Redig, F. and Saada, E.:
  The Abelian sandpile model on an infinite tree.
  \AOPsh\ {\bf 30}, 2081--2107, (2002).

\bibitem{mrs2} Maes, C., Redig, F. and Saada, E.:
  The infinite volume limit of dissipative Abelian sandpiles.
  \CMPsh\ {\bf 244}, 395--417 (2004).

\bibitem{MRSvM} Maes, C., Redig F., Saada, E.~and Van Moffaert, A.:
  On the thermodynamic limit for a one-dimensional sandpile process.
  Markov Process.~Related Fields {\bf 6}  1--22 (2000).

\bibitem{MRS05} Maes, C., Redig, F. and Saada, E.:
  Abelian sandpile models in infinite volume.
  Sankhya, The Indian Journal of statistics {\bf 67}, 634--661 (2005).

\bibitem{MD92} Majumdar, S.N.~and Dhar, D.:
  Equivalence between the Abelian sandpile model and the
  $q \to 0$ limit of the Potts model.
  Physica A {\bf 185} 129--145 (1992).

\bibitem{mrz} Meester, R., Redig, F.~and Znamenski, D.:
  The Abelian sandpile: a mathematical introduction.
  Markov Process.~Related Fields {\bf 6}, 1-22, (2000).

\bibitem{pmnt91} Pemantle, R.:
  Choosing a spanning tree for the integer lattice uniformly.
  \AOPsh\ {\bf 19}  1559--1574, (1991).

\bibitem{priezzh94} Priezzhev, V.B.:
  Structure of two-dimensional sandpile. I. Height probabilities.
  \JSPsh\ {\bf 74}, 955--979 (1994).

\bibitem{priezzh} Priezzhev, V.B.:
  The upper critical dimension of the Abelian sandpile model.
  \JSPsh\ {\bf 98}, 667-684 (2000).

\bibitem{Redig06} Redig, F.:
  Mathematical aspects of the abelian sandpile model.
  Les Houches, Session LXXXIII 2005, A.\ Bovier, F.\ Dunlop, F.\ den Hollander,
  A.\ van Enter and J. Dalibard (eds.), Elsevier, pp. 657--728, (2006).

\bibitem{rosenblatt} Rosenblatt, M.:
  Transition probability operators.
  Proc. Fifth Berkely Symposium Math. Stat. Prob., {\bf 2}, 473-483, (1967).

\bibitem{Wilson} Wilson, D.B.:
  Generating random spanning trees more quickly than the cover
  time. In \emph{Proceedings of the Twenty-Eighth ACM Symposium
  on the Theory of Computing} 296--303. ACM, New York, (1996).

\end{thebibliography}
\end{document}